\documentclass[12pt,twoside]{article} 
\usepackage{amsmath,amssymb,mathrsfs,graphicx}

\newtheorem{theorem}{Theorem}[section]
\newtheorem{lemma}[theorem]{Lemma}
\newtheorem{proposition}[theorem]{Proposition}
\newtheorem{corollary}[theorem]{Corollary}
\newtheorem{definition}[theorem]{Definition}

\newtheorem{rmrk}[theorem]{Remark}

\DeclareMathAlphabet{\mathbfit}{OML}{cmm}{b}{it}

\makeatletter
\@addtoreset{equation}{section}
\makeatother

\newcommand{\fig}[3] {
\medskip\smallskip
\begin{figure}[htb]
  \centering
  \includegraphics[width=#2]{mm-#1.pdf}
  \begin{minipage}[t]{0.80\linewidth} 
    \caption{#3}
    \protect\label{#1}
  \end{minipage}
\end{figure}
\medskip
}

\newenvironment{remark}
{\begin{rmrk} \em}
{\end{rmrk}}

\newcommand{\fn} {function}
\newcommand{\me} {measure}

\newcommand{\erg} {ergodic}
\newcommand{\sy} {system}
\newcommand{\hyp} {hyperbolic}
\newcommand{\pr} {probability}
\newcommand{\dsy} {dynamical system}
\renewcommand{\o} {orbit}

\newcommand{\R} {\mathbb{R}}
\newcommand{\C} {\mathbb{C}}

\newcommand{\Z} {\mathbb{Z}}
\newcommand{\N} {\mathbb{N}}
\newcommand{\qed} {\hfill {\small Q.E.D.} \par\medskip}
\newcommand{\skippar} {\par\medskip}
\newcommand{\ds} {\displaystyle}
\newcommand{\proof} {\noindent \textsc{Proof.} }
\newcommand{\proofof}[1] {\noindent \textsc{Proof of {#1}.} }
\newcommand{\article}[3] {\textsc{{#1}}, {\itshape {#2}}, {{#3}}.}
\newcommand{\book}[3] {\textsc{{#1}}, {\itshape {#2}}, {{#3}}.}
\newcommand{\vol} {\textbf}
\newcommand{\eps} {\varepsilon}
\newcommand{\ph} {\varphi}
\newcommand{\rset}[2] {\left\{ #1 \: \left| \: #2 \right. \! \right\} }
\newcommand{\lset}[2] {\left\{ \left. \! #1 \: \right| \: #2 \right\} }

\newcommand{\symmdiff} {\triangle}

\newcommand{\into} {\longrightarrow}

\renewcommand{\iff} {if and only if\ }
\renewcommand{\emptyset} {\varnothing}

\newcommand{\m} {mixing}
\newcommand{\ob} {observable}

\newcommand{\br} {\tau}   % branch fn
\newcommand{\leb} {m}   % Lebesgue measure
\newcommand{\bj} {\mathbfit{j}}
\newcommand{\bor} {\mathscr{B}}   % Borel sigma-algebra
\newcommand{\sca} {\mathscr{A}}   % generic sigma-algebra
\newcommand{\sci} {\mathscr{I}}   % invariant sigma-algebra
\newcommand{\sct} {\mathscr{T}}   % tail sigma-algebra
\newcommand{\mar} {\mathscr{M}}   % Markov sigma-algebra 
\newcommand{\trn} {\sigma}   % translation operator
\newcommand{\itwo}[2] {\mar^{#1}[#2]}   % element of nth partition containing x
\newcommand{\ione}[1] {\mar[#1]}   % element of original partition containing x
\newcommand{\avg} {\overline{\mu}}   % IV average
\newcommand{\scv} {\mathscr{V}}
\newcommand{\iv} {V \nearrow \R}
\newcommand{\ivlim} {\lim_{\iv}}
\newcommand{\ivlimtwo} {\lim_{{\iv} \atop {n \to \infty}}} 
\newcommand{\go} {\mathcal{G}}   % global observables
\newcommand{\lo} {\mathcal{L}}   % local observables
\newcommand{\bigj} {\kappa}
\newcommand{\bigjj} {\bar{\kappa}}
\newcommand{\con} {\mathcal{C}}   % conservative part
\newcommand{\dis} {\mathcal{D}}   % dissipative part
\newcommand{\sa} {\mathbb{S}_a}   % circle of length a
\newcommand{\aiso} {\aleph_\mathrm{iso}}
\newcommand{\ater} {\aleph_\mathrm{ter}}
\newcommand{\lavg} {\overline{\leb}}   % IV Lebesgue average
\newcommand{\qrw} {\mathcal{Q}}   % transition matrix of rw
\newcommand{\E} {\mathbb{E}}
\newcommand{\inv} {\mathcal{I}}   % invariant hull 
\newcommand{\ic} {{\inv_\con}}
\newcommand{\id} {{\inv_\dis}}
\newcommand{\icd} {{\inv_{\con\dis}}}

\renewcommand{\i} {\imath}

\begin{document}

\title{\textbf{Uniformly expanding Markov maps of the real line:
exactness and infinite mixing}}

\author{\textsc{Marco Lenci}
\thanks{
Dipartimento di Matematica, 
Universit\`a di Bologna, 
Piazza di Porta S.\ Donato 5, 
40126 Bologna, Italy.
E-mail: \texttt{marco.lenci@unibo.it} } 
\thanks{Istituto Nazionale di Fisica Nucleare, 
Sezione di Bologna, Via Irnerio 46,
40126 Bologna, Italy.}
}

\date{Final version for \\
\emph{Discrete and Continuous Dynamical Systems A} \\[10pt]
February 2017}

\maketitle

\begin{abstract}
  We give a fairly complete characterization of the exact components
  of a large class of uniformly expanding Markov maps of $\R$. Using
  this result, for a class of $\Z$-invariant maps and finite
  modifications thereof, we prove certain properties of infinite
  mixing recently introduced by the author.

\end{abstract}

\tableofcontents

\section{Introduction}
\label{sec-intro}

Uniformly expanding Markov maps of the interval represent a paradigm
for chaotic \dsy s. They make a fairly large class of non-trivial
maps, and possess the standard ingredient for chaos, namely, \hyp
ity---insofar as expansivity can be understood as the one-dimensional
version of \hyp ity. On the other hand, they are simple enough to be
more or less fully understood via the techniques of the modern theory
of \dsy s; see the excellent textbook by Boyarsky and G\'ora
\cite{bg}.

In infinite \erg\ theory, the analogues of such maps are the uniformly
expanding Markov maps of $\R$ (see, e.g., Fig.~\ref{pwe} further
down).  Not much is known about them, at least to this author. We are
especially interested in ``translation indifferent'' maps, that is,
maps whose local properties are uniformly bounded throughout $\R$. By
way of counterexample, we are not interested in Boole's transformation
\cite{aw}, which is very close to the identity outside of a compact
set, or in Bugiel's maps \cite{b1, b2}, which are designed to preserve
a finite \me.

In this note we are concerned with the \m\ properties of a very
general class of uniformly expanding Markov maps of the real line.

Initially, we consider the exactness property, which is a strong
notion of mixing that has the advantage of being defined in the same
way in both finite and infinite \erg\ theory. We prove a series of
results that characterize the \erg\ and exact components of a map in
terms of its combinatorics relative the Markov partition. The
characterization is rather precise outside of the invariant set where
the \o s escape to $\pm \infty$. Understandably, the ways in which an
\o\ can escape are many and not easily classifiable. With a few extra
assumptions, however, we are able to give a comprehensive description
of the exact components of this set as well.  A byproduct of all these
results is a number of easily checkable sufficient conditions for the
exactness of a uniformly expanding Markov map.

Later, we apply the notions of mixing for infinite-\me-preserving 
\dsy s (for short, \emph{infinite mixing}) recently introduced by the
author in \cite{limix}. We present these notions, within the present
scope, in Section \ref{sec-im} below and refer the reader to
\cite{limix, liutam, lpmu} for a more thorough discussion. (The last
reference, in particular, uses a more intuitive notation and contains
several results that are used in this article.)

For this part, we specialize to a much narrower but still nontrivial
class of maps. We consider both \emph{quasi-lifts} of expanding circle
maps, i.e., piecewise smooth, translation invariant maps $\R \into
\R$, whose quotient on a fundamental domain is an expanding map of the
circle (see Fig.~\ref{ql}), and \emph{finite modifications} thereof,
namely, maps that coincide with a quasi-lift of an expanding circle
map outside a bounded domain (see Fig.~\ref{fm}). In both cases, we
prove versions of \emph{global-local \m} and \emph{global-global \m}.
Very loosely, global-local \m\ means that any \emph{global \ob}
(roughly, a bounded \fn) and any \emph{local \ob} (an integrable \fn)
decorrelate in time. Global-global \m\ means that the same happens for
any two global \ob s.

Of course, we have stronger results for the more specialized class of
\sy s, that is, the quasi-lifts. In particular, we prove a property,
denoted \textbf{(GLM2)}, that can be recast like this:
For any global \ob\ $F: \R \into \C$ and any probability \me\ $\nu$,
absolutely continuous w.r.t.\ a reference infinite \me\ $\mu$, if
$T_*^n \nu$ denotes the push-forward of $\nu$ via the map $T^n$, then
\begin{equation}
  \label{glm2-intro}
  \lim_{n \to \infty} T_*^n \nu(F) = \avg(F),
\end{equation}
where $\avg(F)$ represents, in a sense that is specified below, the
average of $F$ over $\R$, relative to $\mu$. Thus (\ref{glm2-intro})
can be regarded as a sort of weak convergence of $T_*^n \nu$, the
statistical state of the \sy\ at time $n$, to the ``equilibrium
state'' $\avg$, which is independent of the initial condition $\nu$.
The global \ob s play the role of test \fn s.

We are unable to prove this strong property for all finite
modifications of quasi-lifts of circle maps, but we certainly believe
it to be true for a large class of such \sy s. For this reason, we
give an example for which a very strong version of \textbf{(GLM2)} can
be indeed be proved. In a sense which will be explained below, cf.\
Section \ref{subs-rw}, this example represents a random walk in $\Z$.

\skippar

This is how the paper is organized. In Section \ref{sec-setup} we
introduce our maps and present several results on their exact
components, from the more general statements to the ones that require
extra assumptions. In Section \ref{sec-es} we focus on three
subclasses of maps: the quasi-lifts of expanding circle maps, their
finite modifications and the random walks. In Section \ref{sec-im} we
give our definitions of infinite \m\ and apply them to the \sy s of
Section \ref{sec-es}. The proofs of all the main results are found in
Section \ref{sec-proofs}. The Appendix comprises two sections: in the
first we discuss the importance of some of our assumptions and in the
second we place a standard distortion argument used in Section
\ref{sec-proofs}.

\bigskip\noindent
\textbf{Acknowledgments.}\ I thank Stefano Isola and Jooyoun Hong for
helpful discussions, Sara Munday for her help during the preparation
of the paper, and an anonymous referee for prodding me to make big and
small improvements throughout the paper.  I also acknowledge the
hospitality of the Courant Institute of Mathematical Sciences at New
York University, where part of this work was done. This research is
part of my activity within the \emph{Gruppo Nazionale di Fisica
Matematica} (INdAM, Italy). It was also partially supported by PRIN
Grant 2012AZS52J\underline{ }001 (MIUR, Italy).

\section{Setup and exactness}
\label{sec-setup}

In this section we give the precise definition of our maps of interest
and give a number or results about their exactness properties.  We
start with a characterization of the exact components which intersect
the conservative part of the phase space. Their complement, the
`escape part', has a more complicated dynamics: the results for this
part will require more assumptions, cf.\ (A5)-(A7) in Section
\ref{subs-ex}.

\skippar

Let $( a_j )_{j \in \Z}$ be a collection of real numbers such that
$\lim_{j \to \pm\infty} a_j = \pm\infty$ and
\begin{itemize}
\item[(A1)] $\exists \theta > 0$ such that $0 < a_{j+1}-a_j \le 
  \theta$, $\forall j \in \Z$.
\end{itemize}
Let $I_j := [a_j, a_{j+1}]$. We call $\{ I_j \}_{j \in \Z}$ a
partition of $\R$ even though formally it is not---the substance of
what we discuss in this paper would not change if we made the cleaner
yet more cumbersome choice $I_j := [a_j, a_{j+1})$. Let us denote by
$\mar$ the $\sigma$-algebra generated by the $I_j$.

We consider $T: \R \into \R$, a surjective Markov map relative to $\{
I_j \}$. More precisely,
\begin{itemize}
\item[(A2)] $T |_{(a_j, a_{j+1})}$ has a unique extension $\br_j: I_j
  \into J_j$, which is twice differentiable and bijective onto $J_j
  \in \mar$. Equivalently, $J_j := \bigsqcup_{k=k_{1j}}^{k_{2j}} I_k$,
  for some $k_{1j} \le k_{2j}$.
\end{itemize}
Notice that the above implies that $T$ is two-sided non-singular
w.r.t.\ the Lebesgue \me\ $\leb$. This means that, for all Borel sets
$A$, $\leb(T^{-1} A) = 0 \Leftrightarrow \leb(A) = 0$.

Let $\br_j'$ and $\br_j''$ denote, respectively, the first and second
derivatives of $\br_j$. Then:
\begin{itemize}
\item[(A3)] $\exists \lambda > 1$ such that $| \br_j' | \ge \lambda$,
  $\forall j \in \Z$;
\item[(A4)] $\exists \eta > 0$ such that $\ds \left| \frac{\br_j''} 
  { (\br_j')^2} \right| \le \eta$, $\forall j \in \Z$.
\end{itemize}

An example of a map satisfying (A1)-(A4) is shown in 
Fig.~\ref{pwe}.

\fig{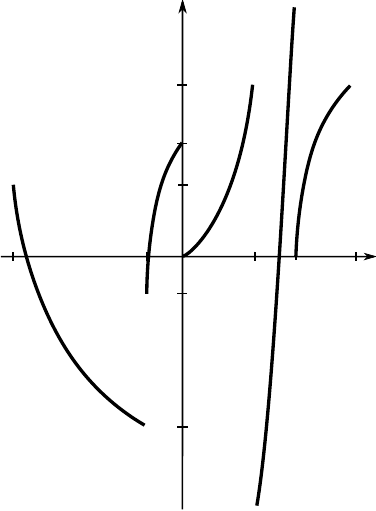}{6cm}{A uniformly expanding Markov map $\R \into \R$.}

\subsection{Classification of Markov intervals}
\label{subs-class}

Throughout the paper we use the following 

\skippar\noindent
\textbf{Convention.}\ All
equalities/inclusions of subsets of $\R$ are intended mod
$\leb$ within $\bor$, the Borel $\sigma$-algebra of $\R$.  In
particular, the strict inclusion $A \subset B$ means $\leb( A \cap B^c
) = 0$ and $\leb( B \setminus A ) > 0$.

\skippar

For a Markov map, the elements of its Markov partition---in our case,
the intervals $I_j$---can be classified in analogy with the states of
a Markov chain; cf., e.g., \cite[Chap.\ VIII]{s} or \cite[Chap.\
4]{g}.  In recalling the definitions below, we will say indifferently
that the interval $I_j$ possesses a certain property or that the state
$j$ possesses that property.

The \emph{transition matrix associated to $T$} is the stochastic
matrix $\mathcal{P} = \mathcal{P}_T := (p_{jk})_{j,k \in \Z}$, where
$p_{jk} := \leb( T^{-1} I_k \,|\, I_j)$. The surjectivity of $T$
implies that, for every $k$, there exists $j$ such that $p_{jk} > 0$.
We denote by $p_{jk}^{(n)}$ the entries of $\mathcal{P}^n$. The Markov
property (A2) implies that
\begin{equation} \label{p-jk-n}
  p_{jk}^{(n)} > 0 \ \Longleftrightarrow \ \leb( T^{-n} I_k \,|\, I_j)
  > 0 \ \Longleftrightarrow \ T^n I_j \supset I_k. 
\end{equation}
When the above occurs for some $n \in \Z^+$, we say that the interval
$I_k$ is \emph{accessible} from $I_j$, or that $I_j$ \emph{feeds}
$I_k$.
  
The intervals $I_j, I_k$ (or the states $j, k$) are called
\emph{communicating} if each one is accessible form the other.  By
convention, we declare that $I_j$ communicates with itself.  This
establishes an equivalence relation on $\Z$. The corresponding
equivalence classes are referred to as the \emph{(communicating)
  classes} of $T$, and are denoted $\Z_\alpha$, with $\alpha \in
\aleph$, some countable set.  We also call
\begin{equation} \label{m-alpha}
  M_\alpha := \bigsqcup_{j \in \Z_\alpha} I_j
\end{equation}
the \emph{set associated to $\Z_\alpha$}. If $\Z$ is one (hence the
only) communicating class, $T$ is called \emph{irreducible}.
    
If $I_j$ feeds some $I_k$, but the viceversa does not hold, we say
that $I_j$ is \emph{inessential}. It is easy to see that the property
of feeding a state, or being accessible from a state, carries over
within a communicating class. So we say, for example, that a certain
class is accessible from $I_j$, or it is inessential, etc. The states
or classes that are not inessential are obviously called
\emph{essential}.
   
An essential class $\Z_\alpha$ that is not accessible from any
external states, that is, such that $j \in \Z_\alpha$ and $k \not\in
\Z_\alpha$ imply $p_{jk}^{(n)} = p_{kj}^{(n_1)} = 0$, $\forall n,n_1
\in \Z^+$, is called \emph{isolated}. The index set of all isolated
classes is denoted $\aiso$.  An essential class that is fed by at
least one external state is called \emph{terminal}. For example,
$\Z_\alpha$ can be a terminal class of the state $k$, or of the
inessential class $\Z_\beta$, etc. Notice that an inessential state
can have more than one terminal class, or none---the latter
possibility can occur because the set of states is infinite. The index
set of all terminal classes is denoted $\ater$.

The integer
\begin{equation} \label{d-j}
  d_j := \mathrm{g.c.d.} \rset{n \in \Z^+} {p_{jj}^{(n)} > 0}
\end{equation}
is called the \emph{period} of $I_j$ (if the r.h.s.\ of (\ref{d-j}) is
empty, set $d_j := 0$). Since this definition is the same as for the
Markov chain generated by $\mathcal{P}$, we know \cite[Thm.~4.2.2]{g}
that two intervals in the same class have the same period. This will
be henceforth called the period of the class $\Z_\alpha$, denoted
$d_\alpha$. We say that the \emph{period of $T$} is $d$ if $d_j = d$,
$\forall j \in \Z$. We say that $T$ is \emph{aperiodic} if it has
period 1.
 
Finally, let us endow $\Z$ with a graph structure by declaring that an
edge exists between $j$ and $k$ \iff $p_{jk} + p_{kj} > 0$, that is,
if $T I_j \supset I_k$ or $T I_k \supset I_j$. This is not the usual
graphical representation of the transition probabilities, which is a
directed graph: it is its undirected version.  It is easy to see that
a communicating class of $T$ and all the states feeding it are all
contained in one connected component of this graph. In fact, a
connected component may contain more that one terminal class, but only
one isolated class (which, in that case, coincides with the connected
component). If there is only one connected component, we say that $T$
is \emph{Markov-indecomposable}. Obviously, an irreducible $T$ is
Markov-indecomposable.

\subsection{Exactness properties}
\label{subs-ex}

We denote by $\con$ and $\dis$ the conservative and dissipative parts
of $T$, respectively. It is known that $T^{-1} \con \supseteq \con$
and $T^{-1} \dis \subseteq \dis$ (mod $\leb$, which is implicit by our
convention) \cite[Chap.\ 1]{a}. As is customary in the field of
non-singular \dsy s, a set $A$ is called \emph{invariant} relative to
$T$, or \emph{$T$-invariant}, if $T^{-1} A = A$. The set
\begin{equation} \label{inv-hull}
  \inv(A) = \inv_T(A) := \bigcup_{k \in \Z^+} \bigcup_{k \in \N} 
  T^{-k} T^n A
\end{equation}
is called the \emph{invariant hull} of $A$ w.r.t.\ $T$. It is the
smallest $T$-invariant set containing $A$. Let us define the following
invariant sets:
\begin{align} 
  \icd &:= \inv(\con) \cap \inv(\dis); \label{icd} \\
  \ic &:= \inv(\con) \setminus \inv_{\con \dis};  \label{ic} \\ 
  \id &:= \inv(\dis) \setminus \inv_{\con \dis};  \label{id} 
\end{align}
Since $\con$ and $\dis$ are complementary, $\R = \ic \sqcup \id \sqcup
\icd$. We may call $\ic$ the conservative-invariant part, $\id$ the
dissipative-invariant part, and $\icd$ the mixed part of $\R$. We will
see below how $\id$ also deserves the name of `escape part'.

\begin{definition} \label{def-ec} 
  In the present context, a set $A \subseteq \R$ is called an
  \emph{\erg\ component} of $T$ if $\leb(A)>0$, $A$ is $T$-invariant
  and $A$ has no $T$-invariant subset of strictly smaller \me. It is
  called an \emph{exact component} of $T$ if it is an \erg\ component
  and $T|_A$ is exact.
\end{definition}  

Observe that the above is a rather stringent definition of \erg\
component, not allowing for zero-\me\ \erg\ components, which can be
defined via the Ergodic Decomposition Theorem \cite[\S2,2]{a}. On the
other hand, the next proposition shows that the union of all null
invariant sets of $\inv(\con) = \ic \sqcup \icd$ is negligible.

Given $x \in \R$, recall the definition of $\omega(x)$, the
$\omega$-limit set of $x$: it is the set of all the accumulation points
of $( T^n(x) )_{n \in \N}$ \cite[Chap.\ 5]{w}. Call $\Omega$ the set
of all $x \in \R$ with a non-empty $\omega(x)$. Evidently, $\Omega$
is measurable and invariant.

\begin{proposition} \label{prop-omega}
  $\Omega = \inv(\con) = \ic \sqcup \icd$. Furthermore, $\Omega$ is
  decomposed mod $\leb$ into a countable number of (positive-\me)
  \erg\ components.
\end{proposition}

\proof Section \ref{sec-proofs}.
\bigskip

The following results concern the exactness properties of
$T|_\ic$ and $T|_\icd$.

\begin{theorem} \label{thm-c} 
  Under the assumptions (A1)-(A4), $\ic$ is made up of at most
  countably many \erg\ components of $T$, denoted $E_\alpha$.  (Here
  $\alpha$ is just a generic index; see however Proposition
  \ref{prop-ic-c}.)  The periods of all $I_j \subset E_\alpha$ are the
  same: we denote them $d_\alpha$. Also, $E_\alpha$ splits into
  $d_\alpha$ exact components of $T^{d_\alpha}$, denoted $E_{\alpha,
  i}$, with $0 \le i \le d_\alpha - 1$. Each $E_{\alpha, i} \in \mar$.  These 
  are `cyclic components' of $T$ in that $T E_{\alpha, i} = E_{\alpha, i+1}$, 
  for $i \in \{ 0, 1, \ldots, d_\alpha-2 \}$, and $T E_{\alpha, d_\alpha-1} = 
  E_{\alpha,0}$.
\end{theorem}

\proof Section \ref{sec-proofs}.

\begin{corollary} \label{cor-c} 
  The following holds:
  \begin{itemize}
  \item[(a)] $\ic \in \mar$.
  \item[(b)] Every $E_\alpha \in \mar$.
  \item[(c)] Every $I_j \subset \ic$ belongs to an $\mar$-measurable
    exact component of $T^{d_j}$.
  \item[(d)] If $T$ is Markov-indecomposable and aperiodic with a
    non-null $\ic$, then it is conservative, irreducible and exact.
  \end{itemize}
\end{corollary}

\proof Assertion \emph{(b)} follows trivially from Theorem
\ref{thm-c}.  So does \emph{(a)} from \emph{(b)}. Assertion \emph{(d)}
is also easy: if $\leb(\ic) > 0$ then $\ic$ contains at least one
\erg\ component $E_\alpha$, which contains at least one interval
$I_j$.  The Markov-indecomposability of $T$ implies that $E_\alpha =
\inv_T(I_j)$ intersects $I_k$, $\forall k \in \Z$. But $E_\alpha \in
\mar$, whence $E_\alpha = \R$, which cannot be split in smaller cyclic
components by aperiodicity. Therefore, $T$ is irreducible and exact,
and $\R = \con$. Finally, not only does \emph{(c)} follow from the
theorem, as is apparent, but the viceversa holds as well.  This will
be shown in the proof of Theorem \ref{thm-c} in Section
\ref{sec-proofs}.  
\qed

As intuition suggests, the \erg\ components of $T$ have much to do
with the communicating classes introduced in Section
\ref{subs-class}. In the remainder of this section we will establish
relations between the two. For the moment, let us remark that the set
of all states that either belong to or feed an essential class
does not necessarily equal $\Z$, for there might be inessential states
which have no terminal class. The collection of all the latter states
will be denoted $\Z_\infty$.

The next two propositions assume (A1)-(A4) and use the notation of
Section \ref{subs-class}. They will be proved in Section
\ref{sec-proofs}.

\begin{proposition} \label{prop-ic-c} 
  Each \erg\ component $E_\alpha \subseteq \con$ equals $M_\alpha$,
  cf.\ (\ref{m-alpha}), for some $\alpha \in \aiso$.  (Hence the
  integer $d_\alpha$ of Theorem \ref{thm-c} is the period of
  $\Z_\alpha$.)  Viceversa, if $\alpha \in \aiso$ and $\# \Z_\alpha <
  \infty$, then $M_\alpha$ is an \erg\ component $E_\alpha \subseteq
  \con$. If $\alpha \in \aiso$ and $\# \Z_\alpha = \infty$, then
  $M_\alpha$ is either an \erg\ component $E_\alpha \subseteq \con$,
  or a $T$-invariant subset of $\id$.
\end{proposition}

For $\alpha \in \ater$, set
\begin{align}
  & E_\alpha := \bigcup_{n\in\N} T^{-n} M_\alpha; 
  \label{d-alpha} \\
  & T_\alpha := T |_{M_\alpha}: M_\alpha \into M_\alpha.
  \label{t-alpha}
\end{align}
Observe that the definition (\ref{d-alpha}) is consistent with the
statements of Proposition \ref{prop-ic-c}, that is, with the case
$\alpha \in \aiso$. In such case, in fact, (\ref{d-alpha}) reduces to
$E_\alpha = M_\alpha$.

\begin{proposition} \label{prop-ic-d}
  Each \erg\ component within $\icd$ is of the form $E_\alpha$, cf.\
  (\ref{d-alpha}), for some $\alpha \in \ater$. Also, $M_\alpha
  \subseteq \con$, $E_\alpha \setminus M_\alpha \subseteq \dis$, and,
  for a.e.\ $x \in E_\alpha$, $\omega(x) = M_\alpha$.  The map
  $T_\alpha$ defined in (\ref{t-alpha}) is conservative and \erg, and
  $M_\alpha$ splits into $d_\alpha$ exact components of
  $T_\alpha^{d_\alpha}$, which are cyclic in the sense of Theorem
  \ref{thm-c} (again, $d_\alpha$ is the period of $\Z_\alpha$).
  
  Viceversa, if $\alpha \in \ater$ and $\#\Z_\alpha < \infty$, then
  $E_\alpha$ is an \erg\ component contained in $\icd$, with the above
  properties.  If $\alpha \in \ater$ and $\# \Z_\alpha = \infty$, then
  $E_\alpha$ is either an \erg\ component contained in $\icd$, with
  the above properties, or a $T$-invariant subset of $\id$.
\end{proposition}

\begin{corollary} \label{cor-ic-d}
  $\con \in \mar$.
\end{corollary}

\proof Propositions \ref{prop-ic-d} and \ref{prop-omega} show that
$\icd \cap \con$ is the union of a countable number of $M_\alpha \in
\mar$, with $\alpha$ ranging in a subset of $\ater$. But $\con = \ic
\sqcup (\icd \cap \con)$.  Corollary \ref{cor-c}\emph{(a)} completes
the proof.  
\qed

Notice that the `mixed \erg\ components' $E_\alpha \subseteq \icd$
need not be $\mar$-measurable. For example, if $j$ feeds both
$\Z_\alpha$ and $\Z_\beta$, then part of $I_j$ will belong to
$E_\alpha$ and part to $E_\beta$.

Also observe that characterizing the exact components of
$T_\alpha^{d_\alpha}$ on the $\omega$-limit set $M_\alpha$ gives 
complete information about the exactness properties of $T$ on 
$E_\alpha$. In fact, if $M_{\alpha,i}$ ($0 \le i \le d_\alpha-1$) denote 
the exact components of $T_\alpha^{d_\alpha}$ inside $M_\alpha$, 
then
\begin{equation}
  E_{\alpha,i} := \bigcup_{n\in\N} T^{-n d_\alpha} M_{\alpha,i}
\end{equation}
are cyclic sets for $T |_{E_\alpha}$, on each of which the
$(d_\alpha)^\mathrm{th}$ power of the map is exact. This can be seen
as follows. If $A$ is a positive-\me\ subset of $E_{\alpha,i}$, there
exist $B \subseteq A$, $\leb(B) > 0$, and $N \in \N$ such that,
$\forall n \ge N$, $T^{n d_\alpha} B \subseteq M_{\alpha,i}$. For all
such $n$, however, since $T^{d_\alpha}$ is exact on $M_{\alpha,i}$,
\begin{equation}
  \bigcup_{k \in \N} T^{-k d_\alpha} T^{k d_\alpha} T^{n d_\alpha}
  B = M_{\alpha,i},
\end{equation}
whence
\begin{equation} \label{extra-10}
  \bigcup_{j \in \N} T^{-j d_\alpha} T^{j d_\alpha} B = E_{\alpha,i}.
\end{equation}
Since $B \subseteq A$, (\ref{extra-10}) holds as well with $A$ in the
place of $B$. This proves that any positive-\me\ subset of
$E_{\alpha,i}$, in the tail $\sigma$-algebra of $T^{d_\alpha}$, is the
entire $E_{\alpha,i}$.

\skippar

Understandably, if the \dsy\ preserves the Lebesgue or a similar \me,
the mixed part of the reference space, whose dynamics is dissipative
in the basin of attraction of a conservative set, must be null:

\begin{proposition} \label{prop-no-do} 
  If $T$ preserves a \me\ $\mu$ equivalent to $\leb$ (this means, $\mu
  \ll \leb$ and $\leb \ll \mu$), then $\Omega = \con$. Equivalently:
  $\ic = \con$, $\id = \dis$, $\icd = \emptyset$.
\end{proposition}

\proof Pick a \fn\ $\zeta \in L^1(\R, \mu)$ such that $\zeta(x)>0$,
$\forall x \in \R$. By definition of $\Omega$ and Prop.~1.1.6 of
\cite{a},
\begin{equation}
  \Omega \subseteq \rset{x \in \R} {\sum_{n=0}^\infty \zeta \circ
  T^n(x) = \infty} = \con \mbox{ mod } \mu,
\end{equation}
i.e., $\Omega \subseteq \con$ mod $\leb$. The reverse inclusion is
obvious, and in any case implied by Proposition \ref{prop-omega},
which also gives the other claims.  
\qed

We now come to $\id$.

\begin{proposition} \label{prop-ic-dmho}
  $\id \supseteq \bigsqcup_{j \in \Z_\infty} I_j$. 
\end{proposition}

\proof Every $I_j$ with $j \in \Z_\infty$ cannot intersect (hence be
contained in) $\ic$, otherwise, by Proposition \ref{prop-ic-c}, it
would be essential; and cannot intersect $\icd$, otherwise, by
Proposition \ref{prop-ic-d}, it would intersect $E_\alpha$, for some
$\alpha \in \ater,$ implying that $\Z_\alpha$ is a terminal class for
$j$.  
\qed

A thorough description of the \erg\ and exact components of $\id$,
like we have for $\con$ and $\icd$ via Propositions \ref{prop-ic-c}
and \ref{prop-ic-d}, is an arduous task. Counterexamples 1 and 2 of
Section \ref{subs-ce} of the Appendix corroborate this intuition. But
if we agree to a few, simple, extra assumptions on $\id$, the
situation improves a great deal:

\begin{itemize}
\item[(A5)] $\exists \rho>0$ such that, $\forall x \in \id$,
  $|T(x) - x | \le \rho$.
 
\item[(A6)] $\exists \theta_o \in (0, \theta)$ such that, $\forall j
  \in \Z$ with $\leb( I_j \cap \id ) > 0$, $\theta_o \le
  a_{j+1}-a_j$; cf.\ (A1).
   
\item[(A7)] $\forall j \in \Z$ with $\leb( I_j \cap \id ) > 0$,
  $T I_j \supset I_j$.
\end{itemize}

We describe (A5) by saying that $T$ has a bounded action on $\id$. In
view of (A1), (A6) and (A2), this amounts to the existence of
\begin{equation} \label{def-bigj}
  \bigj := \max_{{j \in \Z} \atop {\leb( I_j \cap \id ) > 0}} 
  ( k_{2j} - k_{1j} + 1 ). 
\end{equation}  
In other words, $\bigj \in \Z^+$ is the maximum number of Markov
intervals in $J_j = TI_j$, for all $j$ such that $I_j$ has a
non-negligible intersection with $\id$.

For all such $j$, (A7) assumes in addition that $k_{1j} \le j \le
k_{2j}$.  In this case notice that also $k_{2j} - k_{1j} \ge 1$,
otherwise $T I_j = I_j$, which is prohibited by (A3).

\begin{remark}
  Of course, one does not know the set $\id$ \emph{a priori}.  One can
  however ensure that (A5)-(A7) hold if, for instance, $T$ has a
  bounded action on the whole of $\R$, or on $\dis$; and if the
  conditions of (A6)-(A7) are verified for all $j \in \Z$, or at least
  all $j$ that do not belong to no essential class $\Z_\alpha$, with
  $\# \Z_\alpha < \infty$ (cf.\ Propositions \ref{prop-ic-c} and
  \ref{prop-ic-d}).
\end{remark}  

\begin{theorem} \label{thm-dmho}
  Under the assumptions (A1)-(A5), $\id = \dis_{+\infty} \sqcup
  \dis_{-\infty}$, where
  \begin{displaymath}
    \dis_{\pm\infty} := \rset{x \in \R} {\lim_{n \to \infty} T^n(x) =
    \pm \infty}.
  \end{displaymath}
  Both sets are clearly $T$-invariant. If (A6) holds, then $\leb(
  \dis_{\pm\infty} ) \in \{0, \infty \}$.  If (A7) also holds, then
  each of the two sets is either null or an exact component of $T$.
\end{theorem}

\proof Section \ref{sec-proofs}.

\section{Examples}
\label{sec-es}

We can use the results of Section \ref{subs-ex} to find many examples
of exact maps. We first focus on the simplest cases that are not
piecewise linear. This is done for the purposes of Section
\ref{sec-im}, in which certain properties of infinite \m\ are
verified, for the first time, for truly non-linear maps. After that,
we discuss an important class of piecewise linear maps, the random
walks. They motivate the definitions of Section \ref{subs-class} and
provide examples for the relevance of certain assumptions of Section
\ref{subs-ex}.

\subsection{Quasi-lifts of expanding circle maps}
\label{subs-ql}

Let us consider the case where the elements of the Markov partition
have the same size, e.g., $I_j := [aj, a(j+1)]$, for some $a>0$, and
$T$ acts in the same way on each of them, that is, for $x \in I_j$,
$\br_j(x) = \br_0(x-aj) + aj$. This is equivalent to
\begin{equation} \label{trn-comm}
  T \circ \trn = \trn \circ T,
\end{equation}
where $\trn(x) = \trn_a(x) = x+a$. In other words, $T$ is translation
invariant (by the quantity $a$). If $\br_0: I_0 \into J_0$ is
bijective, twice differentiable with bounded second derivative,
expanding and such that $J_0 = [a k_{1,0} , a (k_{2,0} + 1)]$, 
with $k_{1,0}
\le 0$ and $k_{2,0} \ge 0$---cf.\ (A2)---it is easy to see that all
the conditions (A1)-(A7) are verified. See Fig.~\ref{ql}.

\fig{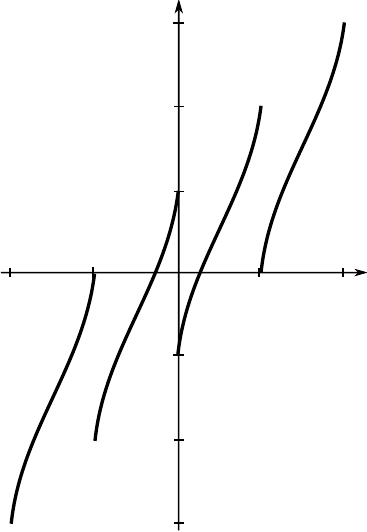}{6cm}{An example of a quasi-lift of an expanding circle map.}

If $\sa$ is the circle constructed by identifying the endpoints of
$[0,a]$, and $T_a: \sa \into \sa$ is defined by $T_a(x) := T(x)$ mod
$a$, we observe that $T_a$ is a uniformly expanding map of the circle,
with bounded distortion and at most one non-regular point, which
happens to be a fixed point.  It thus possesses a number of strong
stochastic properties. In particular, there exists an invariant \me\
$\mu_a$, equivalent to $\leb_a$, the Lebesgue \me\ on $\sa$, which
makes $(\sa, \mu_a, T_a)$ Bernoulli, with exponential decay or
correlations for a large class of \ob s, etc. These results are
proved, e.g., in \cite[Chaps.\ 5 \& 8]{bg}. It may be worth remarking
that $T_a$ is irreducible and aperiodic on $\sa$, so it is exact by
the same arguments proving Theorem \ref{thm-c}.

$T$ is a sort of lift of $\sa$ to $\R$, which we might call
\emph{quasi-lift}. It is apparent that $T$ preserves $\mu$, the
$\sigma$-invariant \me\ whose restriction to $[0,a) \cong \sa$ is
$\mu_a$. Of course, $\mu$ is equivalent to $\leb$. If we set $h_\mu :=
d\mu / d\leb$, the above statements read
\begin{equation} \label{h-inv}
  h_\mu(x) = \sum_{y \in T^{-1}(x)} \frac{h_\mu(y)} {| T'(y) |},
\end{equation}
for all $x \ne aj$ ($j \in \Z$); and $h_\mu(x) = h_\mu(x+a)$,
$h_\mu(x) > 0$, for all $x \in \R$. An important consequence of the
invariance of $\mu$ is that $\con = \ic$ and $\dis = \id$ (Proposition
\ref{prop-no-do}).

A quasi-lift can be thought of as a \emph{$\Z$-extension} of $T_a$,
that is, a self-map $T_\phi$ of $[0,a) \times \Z $ of the form $T_\phi
(y, j) = (T_a(y), j + \phi(y))$, where $\phi: [0,a) \into \Z$
\cite[Chap.~8]{a}. $\phi$ is called (by some) \emph{discrete
displacement}.  In fact, set $\phi (x) := k$ for all $x \in I_0 \cap
T^{-1}(I_k)$ (these sets are intervals and partition $I_0 \cong
[0,a)$, for $k_{1,0} \le k \le k_{2,0}$).  The map $\Psi(y,j) := aj +
y$ defines an isomorphism between the \me\ spaces $( [0,a) \times \Z
\,, \mu_a \otimes \Z )$ and $(\R, \mu)$, and one readily verifies that
$T_\phi = \Psi^{-1} \circ T \circ \Psi$. Therefore, $T_\phi$ preserves
$\mu_a \otimes \Z$.

The quantity 
\begin{equation} \label{def-drift}
  \E_{\mu_a} (\phi) := \frac1 {\mu_a( [0,a) )} \int_0^a \phi \, d\mu_a
\end{equation}
will be called the \emph{drift} of $T$. Notice that there is no harm
in using $\mu$ instead of $\mu_a$ in (\ref{def-drift}). We will do so
throughout the paper.

The assumption (A7) imposes stringent conditions on $\phi$, making the
present \sy s special examples of $\Z$-extensions of expanding circle
maps. Quasi-lifts and similar maps have often been used in nonlinear
physics as toy models for normal and anomalous diffusion; see, e.g.,
\cite{ac1, ac2}, \cite{k, khk}, \cite[Sect.~3.3]{sj}, and references
therein.

\begin{proposition} \label{prop-exa-ql} 
  A quasi-lift of an expanding circle map, as defined above, is exact.
  Furthermore, up to null sets, $\R$ coincides with $\con$,
  $\dis_{+\infty}$, or $\dis_{-\infty}$, depending on the drift
  $\E_\mu (\phi)$ being, respectively, zero, positive, or negative.
\end{proposition}

\proof This result is a corollary of the main theorem of
\cite{ad}. For a much simpler proof, which is self-contained within
the present paper, see Section \ref{sec-proofs}.

\subsection{Finite modifications of quasi-lifts} 
\label{subs-fm}

Slightly more complex examples than the quasi-lifts of circle maps are
the so-called finite modifications of quasi-lifts of circle maps. If
$T_o$ is a quasi-lift as defined earlier, $T: \R \into \R$ is a
\emph{finite modification} of $T_o$ if there exists $k_o \in \Z^+$
such that $T(x) = T_o (x)$, for all $x \not\in
\bigsqcup_{j=-k_o}^{k_o} I_j$. Finite, or \emph{local}, modifications
of translation-invariant \dsy s have been studied in more complicated
contexts as well, such as billiards \cite{laplg, ltyp, dsv}.

We assume that $T$ verifies (A1)-(A7).  Observe that if $T$ is a
finite modification of $T_o$, then $T^n$ is a finite modification of
$T_o^n$, which is a $\Z$-extension by (\ref{trn-comm}). In fact, (A5)
\emph{et seq}.\ show that, if $x \in I_j$, $T(x)$ can land at most
$\bigj-1$ intervals away form $I_j$, hence, for all $x \ne
\bigsqcup_{j=-k_o -n(\bigj-1)}^{k_o +n(\bigj-1)} I_j$, $T^n(x) =
T_o^n(x)$.

Under certain conditions, a finite modification of a quasi-lift is
also exact.

\begin{proposition} \label{prop-exa-fm}
  Let $T$ be a finite modification of a quasi-lift of an expanding
  circle map. It $T$ verifies (A1)-(A7), is Markov-indecomposable and
  preserves a \me\ $\mu$ equivalent to $\leb$, then it is exact.
\end{proposition}

\proof See Section \ref{sec-proofs}.

\begin{remark}
  The significance of the hypotheses of Proposition \ref{prop-exa-fm}
  is clarified by Counterexamples 3 and 4 of Section \ref{subs-ce} of
  the Appendix.
\end{remark}

Under slightly stronger conditions, $T$ verifies a trichotomy similar
to that of Proposition \ref{prop-exa-ql} for quasi-lifts.  The
statement of this result is relatively cumbersome and requires
terminology that is more appropriately introduced later.  For this
reason we present it in Section \ref{sec-proofs}, under the reference
Proposition \ref{prop-exa-fm2}.

We emphasize that the \me\ $\mu$ of Proposition \ref{prop-exa-fm} need
not be the same \me\ preserved by $T_o$, which we henceforth call
$\mu_o$. However, by way of example and in view of later application
(cf.\ Proposition \ref{prop-mix-fm}), we now present a method to
construct finite modifications of quasi-lifts of circle maps which
preserve the original \me.

To simplify things further, we assume that $\mu_o$ is the Lebesgue
\me. The general case can be worked out
in a similar fashion. In view of the notation of Section
\ref{sec-setup}, indicate with $\br_{oj}$ the extension of $T_o
|_{(a_j, a_{j+1})}$ to ${I_j}$, and set $J_{oj} := \br_{oj}(I_j)$.
Define also $\mathcal{J} := \rset{j\in\Z} {J_{oj} \supset I_0}$;
notice that, by (A2), either $J_{oj} \supset I_0$ or $J_{oj} \cap I_0
= \emptyset$, mod $\leb$. Now, pick a $C^2$ \fn\ $\psi: \R \into
\R_0^+$ that is compactly supported in $(a_0, a_1) = (0,a)$. For $j
\in \mathcal{J}$, denote by $\ph_{oj}: J_{oj} \into I_j$ the inverse
\fn\ of $\br_{oj}$, and set $\ph_j := \ph_{oj} + \delta_j \psi$, which
defines a \fn\ on $J_{oj}$. Here $(\delta_j)_{j \in \mathcal{J}}$ is a
collection of numbers so small, in absolute value, that $\ph_j$ is a
monotonic, hence bijective, \fn\ $J_{oj} \into I_j$. (Recall that, by
(A2), $\ph_{oj}$ is monotonic on $J_{oj} \supset I_0$.) Also, they
satisfy
\begin{equation} \label{def-fm10} 
  \sum_{j \in \mathcal{J}} \mathrm{sign} (\ph_{oj}' ) \, \delta_j = 0.
\end{equation}
Finally, define $\br_j := \ph_j^{-1}$, for $j \in \mathcal{J}$, and
$\br_j := \br_{oj}$, otherwise. This determines the map $T$, except at
the points $a_j = aj$, which are negligible.  An example of this
construction is shown in Fig.~\ref{fm}.

\fig{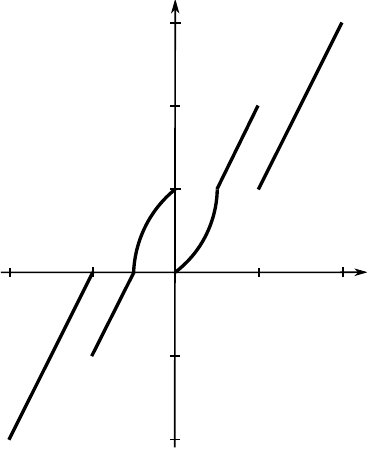}{6cm}{A finite modification of a quasi-lift of a circle map,
  constructed with the procedure given in Section \ref{subs-fm}, for
  the case $\mu_o = \leb$.}

$T$ verifies (A1)-(A2) and (A5)-(A7) by construction. If the
$\delta_j$ are sufficiently small, (A3)-(A4) are verified as well. As
for the invariance of $\mu_o = \leb$, the reader can check that
\begin{equation} 
  \sum_{y \in T^{-1} \{x\}} \frac1 {| T'(y) |} = \sum_{y \in T_o^{-1} \{x\}} 
  \frac1 {| T_o'(y) |},
\end{equation}
at least for all $x \not\in a\Z$. This means that $T$ preserves the 
Lebesgue \me\ \iff $T_o$ does, which was assumed.

\subsection{Random walks}
\label{subs-rw}

A very special family of uniformly expanding Markov maps $\R \into \R$
is given by those representing random walks in $\Z$.

Let $\qrw = (q_{jk})_{j,k \in \Z}$ be the transition matrix of a
random walk in $\Z$, namely, $q_{jk} \in [0,1]$ is the \pr\ that the
walker jumps from the site $j$ to the site $k$.  In line with (A5), we
assume that the walk only admits bounded jumps, i.e., there exists
$\bigjj \in \Z^+$ such that $q_{jk} = 0$, for all $|k-j|>\bigjj$
(although the construction we give below can be easily adapted to the
case of unbounded jumps, cf.\ \cite{lrwwe}).

The map $T = T_\qrw$ associated to the above random walk is 
defined as follows. For $j\in\Z$ and $k \in \{ j-\bigjj, j-\bigjj+1, \ldots,
j+\bigjj \}$, set
\begin{equation} \label{rw-10}
  I_{jk} := \left[ j + \sum_{i = j-\bigjj}^{k-1} q_{ji} \,,\, j + 
  \sum_{i = j-\bigjj}^k q_{ji} \right],
\end{equation}
with the understanding that, when $k = j-\bigjj$, the first of the
above sums is zero. So $\{ I_{jk} \}_{k=j-\bigjj}^{j+\bigjj}$ is a
partition of $[j, j+1]$ into intervals of length, respectively, $\{
q_{jk} \}_{k=j-\bigjj}^{j+\bigjj}$. Notice that, if $q_{jk}=0$,
$I_{jk}$ reduces to a point. We exclude such degenerate intervals. The
complete collection $\{ I_{jk} \}_{j,k \in \Z}$ is the Markov
partition for our map.  For $x$ in the interior of $I_{jk}$ define
\begin{equation} \label{rw-20} 
  T(x) := \frac1 {q_{jk}} \left( x-j- \sum_{i = j-\bigjj}^{k-1} 
  q_{ji} \right) + k.
\end{equation}
For all other $x$, the definition of $T(x)$ is irrelevant. In other
words, $T$ maps $I_{jk}$ affinely onto $[k, k+1]$, see Fig.~\ref{rw}.

\fig{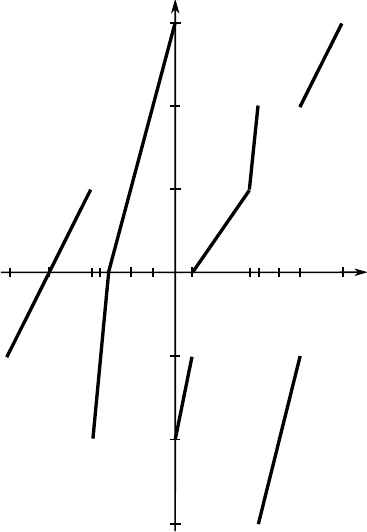}{6cm}{A map $T$ associated a random walk. The marks on the
  abscissa indicate the Markov intervals $I_{jk}$, while those on the
  ordinate represent the intervals $[k, k+1]$.}

A little thinking shows that, if we take a uniformly random $x \in
(k_0, k_0+1) \setminus \bigcup_{n \in \N} T^{-n} \Z$ and look at its
itinerary w.r.t.\ the partition $\{ [k, k+1] \}_{k\in \Z}$, calling
$k_n = k_n(x)$ the unique integer such that $T^n(x) \in (k_n, k_n+1)$,
then $(k_n)_{n \in \N}$ is the random walk on $\Z$ determined by the
transition matrix $\qrw$ and the initial state $k_0$.

Let us observe that $\qrw$ is not the same as $\mathcal{P}_T$, the
transition matrix associated to the map $T$, as presented in Section
\ref{subs-class}. They are however closely related. Denoting by
$\delta_{jk}$ the Kronecker delta, and using that $I_{jk} = [j, j+1]
\cap T^{-1} [k, k+1]$, we get
\begin{equation} \label{rw-50}
\begin{split}
  p_{jk \,,\, j_1 k_1} &:= \leb( T^{-1} I_{j_1 k_1} \,|\, I_{jk} ) \\
  &= \delta_{j_1 , k} \, \frac {\leb ( [j, j+1] \cap T^{-1} [k, k+1] 
  \cap T^{-2} [k_1, k_1+1] )} {\leb ( [j, j+1] \cap T^{-1} [k, k+1] )} 
  \\
  &= \delta_{j_1 , k} \, q_{k,k_1} .
\end{split}
\end{equation}

The map $T$ always verifies the assumptions (A1), (A2), (A4) and
(A5). If $\sup_{j,k} q_{jk} < 1$, (A3) and (A6) are also verified.

On the other hand, (A7) cannot hold in great generality. In fact, it
is verified only if it amounts to the null condition.

\begin{proposition} \label{prop-a7}
  For a map $T$ associated to a random walk such that $\sup_{j,k} 
  q_{jk} < 1$, (A7) can only hold if $\id$ is null.
\end{proposition}

\proof The Markov partition of $T$ is $\{ I_{jk} \}_{j,k}$.  Suppose
by absurd that $\leb (\id) > 0$. There exist $j,k \in \Z$ such that
$\leb( I_{jk} \cap \id ) > 0$. If (A7) holds, $T I_{jk} \supset
I_{jk}$. However, by construction, $T I_{jk} \cap I_{jk} = \emptyset$
(mod $\leb$), $\forall j \ne k$, whence $j=k$.  Since $\id =
\dis_{+\infty} \sqcup \dis_{-\infty}$ and both sets are invariant
(Theorem \ref{thm-dmho}), the trajectories of all $x \in I_{jj} \cap
\id$ must eventually leave $I_{jj}$, implying that $\leb( I_{jk} \cap
\id ) > 0$, for some $k \ne j$. This contradicts what we have just
shown.  
\qed

\begin{remark} \label{rk-rw} 
  This proposition does not imply that---say---random walks with a
  non-zero drift cannot be represented by maps satisfying (A7).  They
  can, only not w.r.t.\ the Markov partition $\{ I_{jk} \}_{j,k}$.
  For example, consider the quasi-lift determined by $\br_0(x) := 3x$,
  where $\br_0$ is the branch of $T$ defined on $[0,1] =: I_0$, cf.\
  Section \ref{subs-ql}. This map represents the homogeneous random
  walk $q_{j,j} = q_{j,j+1} = q_{j,j+2} = 1/3$, $\forall j \in \Z$.
  Clearly $\R = \id = \dis_{+\infty}$. Nonetheless, as discussed in
  Section \ref{subs-ql}, $T$ verifies all (A1)-(A7), relative to the
  Markov partition $\{ I_j = [j, j+1] \}_{j \in \Z}$.
\end{remark}

The following simple result will be useful in the remainder.

\begin{proposition} \label{prop-ds-im}
  The map $T$ associated to the random walk determined by $\qrw$
  preserves the Lebesgue \me\ $\leb$ \iff $\qrw$ is doubly stochastic,
  i.e., $\sum_{j \in \Z} q_{jk} = 1$, $\forall k \in \Z$.
\end{proposition}

\proof We prove this simple proposition by means of the
Perron-Frobenius operator $P = P_T$, which is the operator $L^1(\R,
\leb) \into L^1(\R, \leb)$ uniquely determined by the identity
\begin{equation} \label{pf-10}
  \int_\R (F \circ T) g \, d\leb =  \int_\R F (P g) \, d\leb,
\end{equation}
with $F \in L^\infty(\R, \leb)$ and $g \in L^1(\R, \leb)$.  It is well
known \cite{bg} that, for a.e.\ $x \in \R$,
\begin{equation} \label{pf-20}
  (P g) (x) =  \sum_{y \in T^{-1} \{x\}} \frac{g(y)} {| T'(y) |}.
\end{equation}
For $T$ as in the statement of the proposition, this reads: for all $k
\in \Z$ and $x \in (k,k+1)$,
\begin{equation} \label{pf-30}
  (P g) (x) = \sum_{{j\in\Z} \atop {q_{jk}>0}} q_{jk} \, 
  g (\br_{jk}^{-1}(x)) = \sum_{j\in\Z} q_{jk} \, g (\br_{jk}^{-1}(x)),
\end{equation}
where, in accordance with the notation of (A2), $\br_{jk}$ is the
branch of $T$ defined on $I_{jk}$, cf.\ (\ref{rw-20}).

If we allow (\ref{pf-20})-(\ref{pf-30}) to act on $g \in L^\infty$ as
well, it is clear that $T$ preserves $\leb$ \iff $P1 = 1$, with $1(x)
\equiv 1$ (see also (\ref{h-inv})); that is, \iff $\sum_j q_{jk} = 1$
for all $k$.  
\qed

Markov maps representing random walks are also useful in this paper
for they provide examples which clarify the importance of some of our
earlier assumptions. The reader is referred to Section \ref{subs-ce}
of the Appendix.

\section{Infinite mixing}
\label{sec-im}

In this section we consider the notions of infinite \m\ introduced in
\cite{limix} and further developed in \cite{lpmu}.  We first formalize
them for the case of uniformly expanding maps of the real line and
then apply them to the examples of Sections \ref{subs-ql} and
\ref{subs-fm}.

\subsection{Generalities}
\label{subs-im1}

Consider a ``translation-indifferent'' $T: \R \into \R$. With this
imprecise term we mean that the relevant properties of $T$---e.g.,
expansivity, distortion---are uniform throughout $\R$. In other words,
the map does not single out any special region of $\R$. In this vague
sense, all the examples of Section \ref{sec-es} are
translation-indifferent, even the finite modifications of quasi-lifts,
because the modification does not alter the nature of the map
there. Suppose that $T$ preserves a Lebesgue-absolutely continuous
\me\ $\mu$, which we assume infinite due to translation-indifference.

We call \emph{global \ob} any complex-valued \fn\ $F \in L^\infty(\R,
\mu)$ such that
\begin{equation} \label{im-10}
  \avg(F) := \lim_{r \to \infty} \, \frac1 {\mu([x_0-r, x_0+r])} 
  \int_{x_0 - r}^{x_0 + r} \!\!\!\ F\, d\mu 
\end{equation}
exists uniformly in $x_0$ and independently of it, as the notation
suggests. Clearly, the class of all global \ob s forms a linear space,
containing, for example, the constant \fn s, all \fn s that differ
from a constant by a bounded integrable \fn, or all bounded $F$ with
$\lim_{|x| \to \infty} F(x) < \infty$, etc. Naturally, one is
interested in more complicated \ob s, such as periodic, quasi-periodic
and generally oscillating \fn s: to determine whether they are global
\ob s, one should know $\mu$.

If we restrict $T$ to be a Markov map verifying (A5), we can view this
definition within the general framework presented in \cite{limix,
liutam, lpmu}. We especially refer the reader to \cite{lpmu}, which
uses the same notation as the present paper and contains several
results needed here.

We first assume that $\exists \theta_1, \theta_2 > 0$ such that,
$\forall j \in \Z$,
\begin{equation} \label{im-15}
  \theta_1 \le \mu(I_j) \le \theta_2.
\end{equation}
This makes sense for translation-indifferent \sy s, cf.\ (A1) and
(A6). The collection of sets
\begin{equation} \label{im-20}
  \scv := \lset{ \bigsqcup_{j=k}^\ell I_j } {k \le \ell}.
\end{equation}
is called the \emph{exhaustive family} relative to the Markov
partition of $T$: its elements play the role of ``large boxes'' in
phase space.  Since global \ob s are bounded, it is easy to see that
$F$ verifies (\ref{im-10}) \iff
\begin{equation} \label{im-30}
  \lim_{M \to \infty} \sup_{V \in \scv \atop \mu(V) \ge M} \left| 
  \frac1{\mu(V)} \int_V F \, d\mu - \avg(F) \right| = 0.
\end{equation}
We describe this situation by saying that the average of $F$ over $V
\in \scv$, also denoted $\mu_V (F) := \mu(V)^{-1} \int_V F\, d\mu$,
converges in the \emph{infinite-volume limit} to $\avg(F)$. The
notation
\begin{equation} \label{im-40}
  \ivlim \mu_V (F) = \avg(F)
\end{equation}
is short for (\ref{im-30}). $\avg(F)$ is called the
\emph{infinite-volume average} of $F$.

We also call \emph{local \ob} any complex-valued $g \in L^1(\R, \mu)$.
For any such $g$ we use the customary notation $\mu(g) := \int_\R g\,
d\mu$.

Let us consider two (sub)classes $\go$ and $\lo$ of global and local
\ob s, respectively. Relative to $\go$ and $\lo$, one says that the
\dsy\ $(\R, \mu, T)$ is \m\ of type \textbf{(GLM1)} if, for all $F \in
\go$ and $g \in \lo$, with $\mu(g) = 0$,
\begin{displaymath}
  \lim_{n \to \infty} \mu ((F \circ T^n) g) = 0.
  \eqno{\mathrm{(GLM1)}}
\end{displaymath}
It is \m\ of type \textbf{(GLM2)} if, for all $F \in \go$ and $g \in
\lo$,
\begin{displaymath}
  \lim_{n \to \infty} \mu ((F \circ T^n) g) = \avg(F) \, \mu(g).
  \eqno{\mathrm{(GLM2)}}
\end{displaymath}
It is immediate to see that \textbf{(GLM2)} is equivalent to
(\ref{glm2-intro}) and implies \textbf{(GLM1)}.  As they involve the
pairing of a global and a local \ob, we say that these are two
definitions of \emph{global-local \m}. (There exists another
definition of global-local \m, which is a uniform version of
\textbf{(GLM2)} and is denoted \textbf{(GLM3)} in \cite{lpmu}. We do
not consider it here.)

The following is a trivial consequence of a well-known theorem of Lin
\cite{li}.

\begin{proposition} \label{prop-glm1}
  An exact \dsy\ is {\rm\bfseries (GLM1)}-mixing for any choice of
  $\go \subseteq L^\infty$ and for $\lo = L^1$ (viz.\ any choice of
  $\lo \subseteq L^1$).
\end{proposition}

\proof See \cite[Thm.~3.5\emph{(a)}]{lpmu}.
\bigskip

When we consider the ``decorrelation'' between two global \ob s, we
study the so-called \emph{global-global \m}. We have two definitions
for it.  The \sy\ is called \m\ of type \textbf{(GGM1)} if, for all
$F,G \in \go$, $\avg ( (F \circ T^n) G )$ exists for all sufficiently
large $n$, and
\begin{displaymath}
  \lim_{n \to \infty} \, \avg ( (F \circ T^n) G ) = \avg(F) \, \avg(G).
  \eqno{\mathrm{(GGM1)}}
\end{displaymath}
It is called \m\ of type \textbf{(GGM2)} if, for all $F,G \in \go$,
\begin{displaymath}
  \ivlimtwo \, \mu_V ( (F \circ T^n) G ) = \avg(F) \, \avg(G).
  \eqno{\mathrm{(GGM2)}}
\end{displaymath}
The above limit, which we call \emph{joint infinite-volume and time
limit}, means
\begin{equation} \label{im-60} 
  \lim_{M \to \infty} \sup_{{V \in \scv} \atop {{\mu(V) \ge M} \atop 
  {n \ge M}}} \left| \frac1{\mu(V)} \int_V (F \circ T^n) G \, d\mu - 
  \avg(F) \, \avg(G) \right| = 0.
\end{equation}
The second definition is in essence stronger than the first, as the
following proposition shows.

\begin{proposition} \label{prop-mgg2-1} 
  If $F,G \in \go$ are such that $\avg( (F \circ T^n) G )$ exists for
  all $n$ large enough (depending on $F,G$), then
  \begin{equation} \label{mgg2-1-impl}
    \ivlimtwo \mu_V ((F \circ T^n) G) = b \quad \Longrightarrow 
    \quad \lim_{n \to \infty} \avg ((F \circ T^n) G) = b.
  \end{equation}
  In particular, if the above hypothesis holds $\forall F,G \in \go$,
  then {\rm\bfseries (GGM2)} implies {\rm\bfseries (GGM1)}.
\end{proposition}

\proof See \cite[Prop.~2.3]{lpmu}.
\bigskip

Some of the maps considered in this paper give a good sense of the
relative strength of \textbf{(GGM2)} and \textbf{(GGM1)}, as the
latter property will be trivially verified while the former will
remain an open question; cf.\ Proposition \ref{prop-mix-fm}.

For an in-depth discussion on the meaning and relevance of the above
definitions we refer the reader to \cite{limix}. Here we just point
out that, in order for them to make sense as indicators of
decorrelation, it must be that, for all $F \in \go$ and $n \in \N$,
\begin{equation} \label{im-65}
  \avg(F \circ T^n) = \avg(F). 
\end{equation} 
As shown in \cite{limix}, this is guaranteed by the following
hypothesis: for all $n \in \N$,
\begin{equation} \label{im-70}
  \ivlim \frac{ \mu( T^{-n} V \symmdiff V ) } { \mu(V) } = 0,
\end{equation}
in the sense of the infinite-volume limit, as in (\ref{im-30}).  The
above is easily verified for all \dsy s which verify (A1)-(A2),
(A5)-(A6).  In fact, recalling the definition (\ref{def-bigj}), if $V
= \bigsqcup_{j=k}^\ell I_j$, with $\ell-k$ sufficiently large, it is
easy to see that
\begin{equation} \label{im-80}
  \bigsqcup_{j=k + n (\bigj-1)}^{\ell - n(\bigj-1)} \hspace{-12pt} I_j \ 
  \subset \ T^{-n} V \ \subset \hspace{-4pt} \bigsqcup_{j=k - n 
  (\bigj-1)}^{\ell +  n(\bigj-1)}
  \hspace{-12pt} I_j ,
\end{equation}
whence
\begin{equation} \label{im-90}
  T^{-n} V \symmdiff V \ \subset \hspace{-4pt} 
  \bigsqcup_{j=k - n (\bigj-1)}^{k +  n(\bigj-1)} \hspace{-12pt} I_j 
  \hspace{6pt} \sqcup \hspace{-4pt} \bigsqcup_{j=\ell - n 
  (\bigj-1)}^{\ell +  n(\bigj-1)} \hspace{-12pt} I_j.
\end{equation}
Using (\ref{im-90}) and (\ref{im-15}), we see that the numerator of
(\ref{im-70}) is bounded above by $(4n (\bigj-1) + 2) \theta_2$, while
the denominator is bounded below by $(\ell-k+1) \theta_1$. But the
infinite-volume limit here corresponds precisely to the limit $\ell-k
\to +\infty$ (uniformly in $k, \ell$), whence the assertion.

\subsection{Results for quasi-lifts and their finite modifications}
\label{subs-res-ql-fm}

Now, let $T$ be a quasi-lift of an expanding circle maps, as in
Section \ref{subs-ql}. We are going to show that all the definitions
of infinite \m, both global-local and global-global, are verified for
suitable choices of the global \ob s.

If $\psi$ is either a global or a local \ob, and $k \in \Z^+$, set
\begin{equation} \label{def-ak}
  \mathcal{A}_k \psi := \frac1k \sum_{j=0}^{k-1} \psi \circ \trn^j.
\end{equation}
By (\ref{trn-comm}), this operator commutes with the dynamics, namely
$\mathcal{A}_k (\psi \circ T) = (\mathcal{A}_k \psi) \circ T$. Now
define
\begin{align}
  \label{def-go1}
  \go_1 &:= \rset{F \in L^\infty \! } { \exists F_a = F_a \circ \trn 
  : \lim_{k\to\infty} \| \mathcal{A}_k F - F_a \|_\infty = 0 }; \\
  \label{def-go2}  
  \go_2 &:= \overline{ \mathrm{span}_\C \rset{F \in L^\infty \! } 
  { \exists \beta \in \R : F \circ \trn = e^{\i a \beta}F } },
\end{align}
where the bar denotes closure in the $L^\infty$-norm.  In other words,
$\go_1$ is the space of all essentially bounded \fn s whose
$(a\Z)$-average converges uniformly to a periodic \fn\ (of period
$a$); $\go_2$ is the space generated by the quasiperiodic \fn s
w.r.t. $a\Z$. Clearly, $\go_2 \subset \go_1$ (observe that $\go_1$ is
closed). To see that all these \fn s are global \ob s, we need to
verify that every $F \in \go_1$ possesses an infinite-volume average
$\avg(F)$, in the sense of (\ref{im-30}). In this case, $V$ is of the
form $[ak, a(\ell+1)]$, which gives:
\begin{equation} \label{im-100}
  \frac1{\mu(V)} \int_V F\, d\mu = \frac1{ (\ell-k+1) \mu(I_0) } 
  \int_{ak}^{a(\ell+1)} \!\! F \, d\mu = \int_{I_0} \mathcal{A}_{\ell-k+1} 
  F \circ \trn^k \, d\mu_{I_0}
\end{equation}
which, by the hypotheses on $F$, converges to
\begin{equation} \label{im-110}
  \int_{I_0} F_a \, d\mu_{I_0} =: \avg(F), 
\end{equation}
as $\ell-k \to \infty$ (that is, as $\mu(V) \to \infty$, uniformly in
$V \in \scv$).

Examples of elements of $\go_2$ are the \fn s $E_\gamma(x) := e^{\i
\gamma x}$, $\gamma \in \R$.  An example of $F \in \go_1 \setminus
\go_2$ is given by $F := \sum_{j\in\Z} b_j 1_{I_j}$, with $(b_j)_{j
\in \Z}$ a non-periodic sequence such that $\{ b_{2k}, b_{2k+1}\} =
\{ 0, 1 \}$, for all $k \in \Z$.

We have:

\begin{theorem} \label{thm-mix-ql}
  A quasi-lift of an expanding circle map, as defined in Section 
  \ref{subs-ql}, is:
  \begin{itemize}
  \item[(a)] \m\ of type {\rm\bfseries (GLM1)} for any $\go \subset 
    L^\infty$ and $\lo = L^1$;
  \item[(b)] \m\ of type {\rm\bfseries (GLM2)} w.r.t.\ $\go_1$ and
    $L^1$;
  \item[(c)] \m\ of type {\rm\bfseries (GGM1)} and {\rm\bfseries
    (GGM2)} w.r.t.\ $\go_2$.
  \end{itemize}
\end{theorem}

\proof See Section \ref{sec-proofs}.

\bigskip

The results we have for finite modifications of quasi-lifts are less
satisfactory.

\begin{definition} \label{def-avg-avgo}
  Given $\mu, \mu_o$, two $\sigma$-finite, infinite, \me s on $\R$, we
  write $\avg = \avg_o$ when they admit the same global \ob s and
  coincide on them. This means, for all bounded $F: \R \into \C$,
  $\avg(F)$ exists \iff $\avg_o(F)$ does, and they are equal.
\end{definition}

The above situation can occur, for example, when $h_\mu - h_{\mu_o}
\in L^1(\R, \leb)$, where $h_\mu, h_{\mu_o}$ are the densities of
$\mu, \mu_o$, respectively.

\begin{proposition} \label{prop-mix-fm}
  Let $T$ be a finite modification of a quasi-lift $T_o$ which
  verifies (A1)-(A7), and call $\mu_o$ the \me\ preserved by $T_o$
  (cf.\ Section \ref{subs-fm}). If $T$ is Markov-indecomposable and
  preserves a Lebesgue-equivalent \me\ $\mu$ such that $\avg =
  \avg_o$, then $T$ is:
  \begin{itemize}
  \item[(a)] \m\ of type {\rm\bfseries (GLM1)} for any $\go \subset 
    L^\infty$ and $\lo = L^1$;
  \item[(b)] \m\ of type {\rm\bfseries (GGM1)} w.r.t.\ $\go_2$.
  \end{itemize}
\end{proposition}

\proof See Section \ref{sec-proofs}.

\begin{remark}
  The proof of Proposition \ref{prop-mix-fm} will show that, if one
  drops the hypothesis $\avg = \avg_o$, statement \emph{(a)} still
  holds. As to \emph{(b)}, one still has that $\lim_{n \to \infty}
  \avg_o ( (F \circ T^n) G ) = \avg_o (F) \avg_o (G)$, for all $F,G
  \in \go_2$.
\end{remark}

\subsection{Example: finite modification of a homogeneous 
random walk}
\label{subs-fmrw}

The results of the previous section convince one that the mixing
properties of finite modifications of quasi-lifts are harder to prove,
in general, than those of quasi-lifts.  In particular this holds for
the important notion that we have called \textbf{(GLM2)}.  However,
one expects \textbf{(GLM2)} to hold true for a large class of maps. In
this section we present one such case. Even though we pick a specific
example, the technique generalizes easily to other maps of the same
kind.

Let $T = T_\qrw$ be the map associated to the random walk given by
\begin{equation} \label{fmrw-70}
  \qrw := \frac19 \left(
  \begin{array}{ccccccccccccc}
    \ddots\!\! & \!\!\ddots\!\! & \!\!\ddots\!\! & \!\!\ddots\!\! 
    & \!\!\ddots\!\! 
    & & & & & & & & \\
    & 1 & 2 & 3 & 2 & 1 & & & & & & & \\
    & & 1 & 2 & 3 & 2 & 1 & & & & & & \\\hline
    & & & 1 & 2 & 5 & 1 & 0 & & & & & \\
    & & & & 1 & 1 & 5 & 1 & 1 & & & & \\
    & & & & & 0 & 1 & 5 & 2 & 1 & & & \\\hline
    & & & & & & 1 & 2 & 3 & 2 & 1 & & \\
    & & & & & & & 1 & 2 & 3 & 2 & 1 & \\[-3pt]
    & & & & & & & & 
    \!\!\ddots\!\! & \!\!\ddots\!\! & \!\!\ddots\!\! & \!\!\ddots\!\! 
    & \!\!\ddots 
  \end{array}
  \right),
\end{equation}
with the convention that the entries of $\qrw$ are null outside of the
shown diagonal strip. This matrix is doubly stochastic, so $T$
preserves $\leb$ (Proposition \ref{prop-ds-im}). Also, as indicated in
(\ref{fmrw-70}), its rows $(q_{jk})_k$ fail to be translations of one
another only for $j \in \{ -1, 0, 1 \}$, hence $T$ is a finite
modification of a map $T_o$ representing a homogeneous random walk.
More examples of suitable $\qrw$ can be constructed using the ideas of
\cite[App.~A]{lrwyep}.  There is a sizable literature about finite
modification of translation-invariant random walks. Some recent
references include \cite{ps, n, pp, ip}.

We will show that the \dsy\ $(\R, \leb, T)$ verifies a very strong
instance of \textbf{(GLM2)}. Define $\scv' := \lset{[-k, k] \subset
\R} {k \in \Z^+}$. This exhaustive family is smaller than the one we
have introduced in (\ref{im-20}), which in this specific case reads
$\scv = \lset{[k, \ell+1] \subset \R} {k \le \ell \in \Z}$. In a
sense, up to inessential variations, $\scv'$ is the smallest
collection of sets that make sense as an exhaustive family, because it
contains only one increasing sequence of sets that covers the phase
space $\R$. Therefore, the class of \fn s
\begin{equation} \label{def-go0}
  \go' := \rset{F \in L^\infty(\R, \mu)} {\exists \lavg'(F)
  := \lim_{k \to \infty} \, \frac1{2k+1} \int_{-k}^k F\, d\leb}
\end{equation}
is essentially the largest class of global \ob s one can imagine for
the \dsy\ at hand, because $\lavg'$ is the infinite-volume average
w.r.t.\ $\scv'$; cf.\ (\ref{im-30})-(\ref{im-40}).

\begin{remark}
  It is worthwhile to point out that $\go'$ is not simply a larger
  class of global \ob s than $\go_1$ and $\go_2$. By using $\scv'$ in
  lieu of $\scv$ here, we have changed the notion of infinite-volume
  average from $\lavg$, cf.\ (\ref{im-20})-(\ref{im-40}), to $\lavg'$,
  cf.\ (\ref{def-go0}), and therefore extended the very concept of
  global \ob. Notice that, if $\lavg(F)$ exists, $\lavg'(F) =
  \lavg(F)$.
\end{remark}

\begin{proposition} \label{prop-fmrw-mix}
  The map $T$ defined above is irreducible, conservative and exact.
  Also, in addition to the statements of Proposition
  \ref{prop-mix-fm}, it is {\rm\bfseries (GLM2)} relative to the
  exhaustive family $\scv'$, the class of global \ob s $\go'$ and the
  class of local \ob s $L^1$.
\end{proposition}

\proof See Section \ref{sec-proofs}.
\bigskip

Observe that no form of global-global \m\ can hold for $T$
w.r.t. $\go'$.  This is an instance of a general phenomenon that in
\cite[Sect.~3]{limix} we have called \emph{surface effect}. Very
briefly, a counterexample is constructed by choosing, e.g., $F = G =
\Theta$, the Heaviside \fn, which belongs in $\go'$, with
$\lavg'(\Theta) = 1/2$. Since $T$ has a bounded action, in the sense
of (A5), it is clear that $\Theta \circ T^n(x) = \Theta(x)$, for all
$|x| > \rho n$ (in this particular case one can take $\rho = 2$). So,
for all $n \in \N$,
\begin{displaymath}
  \ivlim \leb_V ((\Theta \circ T^n) \Theta) = \lavg'( \Theta^2 ) = 
  \frac12 \ne \lavg'(\Theta)^2, 
\end{displaymath}
contradicting both \textbf{(GGM1)} and \textbf{(GGM2)}. (Recall that,
in the latter, the convergence in $n$ is uniform w.r.t.\ the one in
$V$, and viceversa.)

The above says no more and no less than: $\scv'$ is the wrong
exhaustive family for the global-global \m\ of quasi-lifts, or finite
modifications thereof; cf.\ \cite[Sect.~3]{limix}.  We still expect
\textbf{(GGM1-2)} to hold for general classes of global \ob s,
relative to the exhaustive family $\scv$.

\section{Proofs}
\label{sec-proofs}

In this section we prove all the results stated in the previous
sections (except for the simplest ones, whose proofs have already been
given). We start by laying out the technical tools that will be needed
in the proofs of all the theorems and propositions of Section
\ref{subs-ex}.

For $\bj = (j_0, j_1, \ldots, j_{n-1}) \in \Z^n$, set
\begin{equation} \label{exa-5} 
  I^{(n)}_\bj := I_{j_0} \cap T^{-1} I_{j_1} \cdots \cap T^{-n+1}
  I_{j_{n-1}}.
\end{equation}
By (A1)-(A3), $\{ I^{(n)}_\bj \}_{\bj \in \Z^n}$ is a Markov partition
for $T^n$ (as always, modulo the endpoints of the intervals) and
$\leb(I^{(n)}_\bj) \le \theta \lambda^n$. We call $\mar^n$ the
generated $\sigma$-algebra. For $n \ge 1$, let $\itwo{n}{x}$ denote
the only element of $\{ I^{(n)}_\bj \}$ such that $x \in \itwo{n}{x}$
(in case $x$ belongs to two such elements, being an endpoint of both,
we make the convention that $\itwo{n}{x}$ is the interval on the
right). $\ione{x}$ will be short for $\itwo{1}{x}$.

\begin{lemma} \label{lem-ex1}
  If $x$ is a density point of $A$, with $\leb(A) > 0$, then
  \begin{displaymath}
    \lim_{n \to \infty} \leb(T^n A \,|\, \ione{T^n(x)} ) = 1.
  \end{displaymath}
\end{lemma}

\proof By the hypothesis on $x$,
\begin{equation} \label{exa-10}
  \lim_{n \to \infty} \leb(A \,|\, \itwo{n}{x} ) = 1.
\end{equation}
Setting $B_{x,n} := \itwo{n+1}{x} \setminus A$, (\ref{exa-10}) is
equivalent to $\lim_{n \to \infty} \leb(B_{x,n} \,|\, \itwo{n+1}{x}) =
0$, whence, by Corollary \ref{cor-dist} (Section \ref{subs-dist} of
the Appendix),
\begin{equation} \label{exa-15}
  \lim_{n \to \infty} \leb(T^n B_{x,n} \,|\, T^n \itwo{n+1}{x}) = 0.
\end{equation}
Since $T^n$ is a bijection $\itwo{n+1}{x} \into \ione{T^n(x)}$ and
$\ione{T^n(x)} \setminus T^n A \subseteq T^n B_{x,n}$, (\ref{exa-15})
implies the lemma.  
\qed

\begin{corollary} \label{cor-ex1}
  For $d \in \Z^+$, suppose that $A$ is $T^d$-invariant, with $\leb(A)
  > 0$; $x$ is a density point of $A$; and $y$ an accumulation point
  of $(T^{dn} (x))_{n\in\N}$.  If $y \in (a_j, a_{j+1})$, for some
  $j\in \Z$, then $I_j \subseteq A$ (this means, as always,
  $\mathrm{mod}\ \leb$). If $y = a_j$ and $(T^{dn} (x))_n$ accumulates
  to $y$ from the right (respectively, left), then $I_j \subseteq A$
  (respectively $I_{j-1} \subseteq A$).
\end{corollary}

\proof Evaluate the limit of Lemma \ref{lem-ex1} on any subsequence
$(dn_i)_{i \in \N}$ such that $T^{dn_i}(x) \to y$ from the right/left,
respectively.  
\qed

\bigskip

\proofof{Proposition \ref{prop-omega}} Clearly $\con \subseteq
\Omega$, and so $\inv(\con) \subseteq \Omega$. The first assertion of
the proposition will be proved once we show the reverse inclusion.
This amounts to show that the invariant set $\inv_{\dis \Omega} := \id
\cap \Omega$ is null.

Suppose not. For $j \in \Z$, let $\Omega_j$ be the set of all $x \in
\Omega$ whose $\omega$-limit satisfies at least one of the following
conditions:
\begin{enumerate}
\item $\omega(x) \cap (a_j, a_{j+1})$ is not empty;
\item $a_j \in \omega(x)$ and $( T^n(x) )_{n \in \N}$ accumulates to
  $a_j$ from the right; 
\item $a_{j+1} \in \omega(x)$ and $( T^n(x) )_{n \in \N}$ accumulates 
  to $a_{j+1} $ from the left.
\end{enumerate}
Clearly, $\Omega_j$ is $T$-invariant, and so is $\inv_{\dis \Omega}
\cap \Omega_j$. Since $\inv_{\dis \Omega}$ has positive \me\ and the
\o s of all its elements accumulate somewhere, there exists $j \in \Z$
such that $\leb(\inv_{\dis \Omega} \cap \Omega_j) > 0$. We claim that
the \o\ of a.e.\ $x \in I_j$ returns to $I_j$ infinitely many times in
the future. If not, there would exist $B \subseteq I_j$, with $\leb(B)
> 0$, and $N \in \N$, such that $I_j \cap \bigcup_{n \ge N} T^n B$
is null. One the other hand, the typical $x \in \inv_{\dis \Omega}
\cap \Omega_j$ is a density point of the same set. Applying Corollary
\ref{cor-ex1} (with $d=1$) we get $I_j \subseteq \inv_{\dis \Omega}
\cap \Omega_j$. Thus $B \subseteq \inv_{\dis \Omega} \cap \Omega_j$,
which is absurd because, by construction, no point of $B$ can belong
to $\Omega_j$.

Observe that, in the terminology of Section \ref{subs-class}, we have
just proved that $I_j$ is essential.

Now take any positive-\me\ $A \subseteq I_j$. A.e.\ $x \in A$ is both
a density point of $A$ and a recurrent point to $I_j$. Choose one such
$x$.  The proof of Lemma \ref{lem-ex1} shows that there exists a
return time $n$ such that $T^n (A \cap \itwo{n+1}{x})$ is so large
within $\ione{T^n(x)} = I_j$ to have a non-null intersection with
$A$. Since $T^n$ acts as a two-sided non-singular bijection
$\itwo{n+1}{x} \into \ione{T^n(x)}$, we obtain
\begin{equation}
  \leb(A \cap T^{-n} A) \ge \leb(A \cap \itwo{n+1}{x} \cap T^{-n} A) > 0.
\end{equation}
We have thus proved that $I_j$ cannot contain wandering sets, which is
a contradiction because $I_j \subseteq \inv_{\dis \Omega} \cap
\Omega_j \subseteq \dis$. This concludes the proof of the first
assertion of Proposition \ref{prop-omega}.

For the second assertion it suffices to prove that every invariant $A
\subseteq \Omega$, with $\leb(A)>0$, contains a positive-\me,
invariant subset $B$ which cannot be further decomposed in invariant
subsets of strictly smaller \me. So, consider one such $A$.  The
previous arguments show that, for some $j \in \Z$, $\leb(A \cap
\Omega_j)>0$ and $I_j \subseteq A \cap \Omega_j$. Therefore $B :=
\inv_T (I_j) \subseteq A$ cannot be further decomposed in smaller
invariant subsets, ending the proof of Proposition \ref{prop-omega}.

We add a few remarks. The above conclusion states that $\inv_T(I_j)$
is an \erg\ component of $T$. If we only take the forward images of
$I_j$, we see that $\bigcup_{n \in \N} T^n I_j \subseteq \con$,
because $I_j \subset \con$, as shown earlier, and $T \con \subseteq
\con$. With reference to the definitions of Section
\ref{subs-class}---see in particular (\ref{m-alpha}) and
(\ref{t-alpha})---let $\alpha$ be the unique index in $\aleph$ such
that $j \in \Z_\alpha$. Then $M_\alpha = \inv_{T_\alpha} (I_j) =
\bigcup_{n \in \N} T^n I_j \subseteq \con$, and $T_\alpha$ is
conservative and \erg.  
\qed

\subsection{Exactness}
\label{subs-pf-ex}

Recall that $\bor$ denotes the Borel $\sigma$-algebra of $\R$. Let us
introduce the other $\sigma$-algebras that we are concerned with.  For
$n \in \N$,
\begin{equation} \label{i-n}
  \sci^n := \rset{A \in \bor} {T^{-n}A = A \mbox{ mod } \leb}
\end{equation}
is the \emph{$T^n$-invariant $\sigma$-algebra}. (From now on, as
declared in Section \ref{sec-setup}, we will always imply `mod
$\leb$'.)  Clearly, if $n$ is a multiple of $k$, $\sci^k \subseteq
\sci^n$. $\sci$ will be short for $\sci^1$.  The \emph{tail
$\sigma$-algebra} is defined to be:
\begin{equation} \label{tail}
  \sct := \bigcap_{n=0}^\infty T^{-n} \bor.
\end{equation}
Of course, $\sci^n \subseteq \sct$, for all $n \in \N$.

Given a $\sigma$-algebra $\sca$ and a Borel $B$, we will denote by
$\sca \cap B := \rset{A \cap B} {A \in \sca}$ the \emph{trace} of
$\sca$ in $B$.

At the core of all exactness proofs will be the following
generalization of a criterion by Miernowski and Nogueira \cite{mn}:

\begin{proposition} \label{prop-mn} 
  Consider the \dsy\ $(X, \sca, \nu, S)$, where $(X, \sca, \nu)$ is a
  $\sigma$-finite \me\ space and $S$ a non-singular endomorphism on it
  (i.e., $\nu(A) = 0$ $\Rightarrow$ $\nu(S^{-1} A) = 0$).  Denote by
  $\sci := \rset{A \in \sca} {S^{-1}A = A \mbox{ mod } \nu}$ and $\sct
  := \bigcap_{n=0}^\infty S^{-n} \sca$, respectively, the invariant
  and tail $\sigma$-algebras.  Clearly, $\sci \subseteq \sct$. If,
  $\forall A \in \sct$ with $\nu(A) > 0$, $\exists n = n(A)$ such that
  $\nu( S^{n+1}A \cap S^n A ) > 0$, then $\sci = \sct$.
\end{proposition}

In other words, under the above hypotheses, the non-null \erg\
components of $S$ are also exact components. The proof of Proposition
\ref{prop-mn}, together with a converse statement, can be found in
\cite[Prop.~A.2]{lsimple}.

\bigskip

\proofof{Theorem \ref{thm-c}} We start by proving that, for all $d \ge
1$, the \erg\ components of $T^d$ within $\ic$ (equivalently, the
\erg\ components of $(T|_\ic)^d$) are $\mar$-measurable.  In fact,
given a $T^d$-invariant $A \subseteq \ic$, consider $j \in \Z$ such
that $\leb( A \cap I_j ) > 0$. Since $A$ is in the conservative part
of $T$ (or $T^d$, which is the same) the typical $x \in A \cap I_j$ is
a density point of $A$ and is recurrent to the interior of $I_j$,
w.r.t.\ $T^d$. Corollary \ref{cor-ex1} shows that $I_j \subseteq A$.
Thus, $\sci^d \cap \ic \subseteq \mar \cap \ic$ and the claim is
proved.

Consider an \erg\ component $E_\alpha \subseteq \ic$. Since $E_\alpha
\in \mar$, it contains whole Markov intervals. Set $\Z_o := \rset{j
\in \Z}{I_j \subseteq E_\alpha}$. For all $j \in \Z_o$, $\inv_T(I_j)
= E_\alpha$ and $I_j$ is essential (by conservativity).  So all these
$I_j$ communicate with each other. Therefore $\Z_o$ is a communicating
class $\Z_\alpha$ and $E_\alpha$ is the corresponding $M_\alpha$, cf.\
(\ref{m-alpha}).  In particular, the periods of all intervals $I_j
\subset E_\alpha$ are the same.  Let $d_\alpha \ge 1$ denote their
common value.  For the sake of simplicity, in the remainder of this
proof we write $d$ for $d_\alpha$.

Now fix $A \subseteq E_\alpha$ with $\leb(A)>0$. We are going to show
that $\exists n \in \Z^+$ such that
\begin{equation} \label{exa-40} 
  \leb ( T^{(n+1)d} A \cap T^{nd} A ) > 0.
\end{equation}
Therefore Proposition \ref{prop-mn} can by applied to $T^d
|_{E_\alpha}$, implying that the \erg\ components of $T^d$ within
$E_\alpha$ are exact.

We start with a simple lemma concerning the return times of a state in
a countable-state Markov chain:

\begin{lemma} \label{lem-rj1}
  Denoting $R_j :=  \{n \in \Z^+ \,|\,p_{jj}^{(n)} > 0 \}$, 
  and recalling that $d_j := \mathrm{g.c.d.}(R_j)$,
  there exists $n \in R_j$ such that $n+d_j \in R_j$. 
\end{lemma} 

\proof By (\ref{p-jk-n}), $R_j$ is an additive set, therefore $R_j -
R_j$ is a $\Z$-ideal. By principality, $R_j - R_j = \Z d_j$, therefore
$\exists n,n_1 \in R_j$ such that $d_j = n_1 - n$.  
\qed

Choose $j$ such that $\leb (A \cap I_j) > 0$. By Lemma \ref{lem-rj1},
\begin{align} 
  & T^{n_1} I_j \supset I_j; \label{exa-50} \\
  &T^{n_1+d} I_j \supset I_j, \label{exa-51} 
\end{align}
for some $n_1 \in \Z^+$.  Thus, $T^d \, T^{n_1} I_j \supset I_j$. The
Markov property of $T^d$ implies that $\exists B \subset T^{n_1} I_j$
such that $T^d |_B : B \into I_j$ is a bijection. (Incidentally, $B
\in \mar^d$.)

In view of (\ref{exa-50}), $\exists \delta>0$ such that every $A'
\subset \R$ with $\leb (A' | I_j) > 1 - \delta$ verifies both
\begin{align} 
  & \leb( T^{n_1} A' | I_j ) > \frac12; \label{exa-60} \\[3pt]
  & \leb( T^{n_1} A' | B ) > 1 - \frac1{2D}, \label{exa-61} 
\end{align}
where $D$ is the distortion constant that appears in Corollary
\ref{cor-dist} of Section \ref{subs-dist}. That corollary, together
with the definition of $B$ and (\ref{exa-61}), gives
\begin{equation} \label{exa-70} 
   \leb( T^{n_1 + d} A' | I_j ) > \frac12.
\end{equation}
Therefore, by (\ref{exa-60}) and (\ref{exa-70}),
\begin{equation} \label{exa-80} 
  \leb ( T^{n_1 + d} A' \cap T^{n_1} A' ) > 0.
\end{equation}

At this point, observe that a.a.\ $x \in A \cap I_j$ are both density
points of $A$ and recurrent to the interior of $I_j$, that is,
$\ione{T^n(x)} = I_j$ for infinitely many $n$. Choose any such $x$. By
Lemma \ref{lem-ex1}, there exists a large enough $n_2$ such that $\leb
(T^{n_2} A | I_j) > 1 - \delta$.  This means that (\ref{exa-80}) can
be applied with $T^{n_2} A$ in the place of $A'$. More precisely,
$\leb ( T^{n_1 + n_2 + d} A \cap T^{n_1 + n_2} A ) > 0$.  By the
non-singularity of $T$, this inequality holds as well if $n_1 + n_2$
is replaced with any $nd \ge n_1 + n_2$, which proves (\ref{exa-40}).

So, the \erg\ components of $T^d$ inside $E_\alpha$ are also exact
components of $T^d$. Let us study them. This part of the proof uses a
standard argument from the classification of states for Markov chains
\cite[Sect.~VIII.2]{s}.

Choose $j$ such that $I_j \subset E_\alpha$. Define $\Z'_0 := \{ j \}$
and, for $\ell \ge 1$,
\begin{equation} \label{exa-90} 
  \Z'_\ell := \rset{ k\in\Z } {p_{jk}^{(\ell)} > 0 }.
\end{equation}
In other words, using also the Markov property of $T$, $\bigsqcup_{k
\in \Z'_\ell} I_k = T^\ell I_j$.  Given $k \in \Z'_\ell$, for all
$n$ such that $p_{kj}^{(n)} > 0$ (and there are infinitely many of
them, because all these Markov intervals are contained in a
conservative \erg\ component of $T$) we have $p_{jj}^{(\ell + n)} >
0$, hence $n \equiv -\ell$ (mod $d$) (because $d_j = d$). This implies
that the sets
\begin{equation} \label{exa-92} 
  \Z''_i := \!\!\!\ \bigcup_{{\ell \in \N} \atop 
  {\ell \equiv i \: (\mathrm{mod} \: d)}}  \!\!\!\ \Z'_\ell,
\end{equation}
defined for $i \in \{ 0, 1, \ldots, d-1 \}$, are pairwise disjoint.
It is clear that their union is the communicating class $\Z_\alpha$
that contains $j$, and so $E_\alpha = M_\alpha$, as observed
earlier. For $i \in \{ 0, 1, \ldots, d-1 \}$, set $E_{\alpha,i} :=
\bigsqcup_{k \in \Z''_i} I_k$. By the definitions
(\ref{exa-90})-(\ref{exa-92}),
\begin{equation} \label{exa-100}
  T E_{\alpha,i} = E_{\alpha,i+1 \: (\mathrm{mod}\ d)}, 
\end{equation}
which is one of the assertions of Theorem \ref{thm-c}. 

It remains to prove that each $E_{\alpha,i}$ is an \erg\ component of
$T^d$. First off, since $E_\alpha = \bigsqcup_{i=0}^{d-1}
E_{\alpha,i}$ is $T$-invariant, and by (\ref{exa-100}), $T^{-d}
E_{\alpha,i} = E_{\alpha,i}$.  Suppose by absurd that $E_{\alpha,i}$
could be split in two non-trivial $T^d$-invariant sets $A,B$. By the
two-sided non-singularity of the map, $T^{-i}A, T^{-i}B$ would be
non-trivial $T^d$-invariant subsets of $E_{\alpha,0}$. They would also
belong to $\mar$, as we have shown at the start of this proof. So one
of them, say $T^{-i}A$, must contain $I_j$. However, by definition of
$E_{\alpha,0}$, that is, by definition of $\Z''_0$, $\bigcup_{n \in \N} T^{nd}
I_j = E_{\alpha,0}$, giving that $T^{-i}A = E_{\alpha,0}$, a
contradiction. This ends the proof of Theorem \ref{thm-c}.  
\qed

\bigskip

\proofof{Proposition \ref{prop-ic-c}} The first assertion was all but
shown in the previous proof: a conservative \erg\ component must be of
the form $M_\alpha$, for some $\alpha \in \aleph$.  Moreover,
$\Z_\alpha$ is an isolated class because $T^{-1} M_\alpha =
M_\alpha$. In other words, $\alpha \in \aiso$.

Viceversa, given $\alpha \in \aiso$, consider the $T$-invariant set $A
:= M_\alpha \cap \Omega$: its Lebesgue \me\ can be either positive or
zero.

If $\leb(A)>0$, Corollary \ref{cor-ex1} (with $d=1$) entails that $I_j
\subseteq A$, for some $j \in \Z_\alpha$. By definition of
communicating class, $\forall k \in \Z_\alpha$, $\exists n \ge 1$ such
that $T^n I_j \supset I_k$. This shows that $A \supseteq M_\alpha$,
whence $A = M_\alpha$. It also shows that $I_j$ has infinitely many
Markov returns. The same arguments as in the proof of Proposition
\ref{prop-omega} prove that there are no wandering sets in
$I_j$. Since $j$ is arbitrary, $M_\alpha \subseteq \con$.  But
$M_\alpha$ is $T$-invariant, whence $M_\alpha \subseteq \ic$.

In the case $\leb(A)=0$, $M_\alpha \subseteq \id$ by Proposition
\ref{prop-omega}.

Lastly, observe that $\#\Z_\alpha < \infty$ entails that $M_\alpha$ is
compact, which implies the first of the two cases above.  
\qed

\bigskip

\proofof{Proposition \ref{prop-ic-d}} Let $E_o$ be an \erg\ component
of $\icd$, whence $\leb(E_o)>0$. Recall the definition of $\Omega_j$
from the proof of Proposition \ref{prop-omega}. For at least for one
$j$, the invariant set $E_o \cap \Omega_j$ has positive \me. The same
arguments as in the aforementioned proof show that $I_j$ is essential
and it is contained in $E_o$.  It follows that $\Z_\alpha$, the
communicating class that contains $j$, is essential, and $M_\alpha
\subseteq E_o$. Then $T M_\alpha = M_\alpha$ and $E_\alpha = E_o$,
for both are the \erg\ component containing $I_j$.

\begin{lemma} \label{lem-w}
  A positive-\me\ $W \subseteq M_\alpha$ is a wandering set for $T$
  \iff it is a wandering set for $T_\alpha$, namely, the conservative
  and dissipative parts of $T$ and $T_\alpha$ coincide within
  $M_\alpha$.
\end{lemma}

\proofof{Lemma \ref{lem-w}} Given $W$ as in the statement of the
lemma, set
\begin{align}
  W_1' &:= T^{-1} W \cap M_\alpha = T_\alpha^{-1} W; \\
  W_1'' &:= T^{-1} W \setminus M_\alpha. 
\end{align}
Of course, $T^{-1}W \cap W \ne \emptyset \ \Leftrightarrow \
T_\alpha^{-1} W \cap W \ne \emptyset$. Also, since $T^{-1} M_\alpha
\supseteq M_\alpha$, one has that $T^{-k} W_1'' \cap M_\alpha =
\emptyset$, $\forall k \ge 0$. Applying the same reasoning with
$T_\alpha^{-1} W$ in the place of $W$ and so on, recursively, we
establish that, for all $n \ge 1$,
\begin{equation}
  T^{-n} W \cap W \ne \emptyset \ \Longleftrightarrow \ T_\alpha^{-n} 
  W \cap W \ne \emptyset,
\end{equation}
which was to be proved.
\qed

By definition of $T_\alpha$, cf.\ (\ref{t-alpha}), $\Z_\alpha$ is an
isolated class for $T_\alpha$, which is irreducible.  Proposition
\ref{prop-ic-c}, applied to $T_\alpha$ instead of $T$, proves that
one, and only one, of the following occurs:
\begin{enumerate}
\item $T_\alpha$ is conservative and enjoys all the \erg\ properties
  listed in the statement of Theorem \ref{thm-c}. By Lemma
  \ref{lem-w}, $M_\alpha \subseteq \con$;
\item $T_\alpha$ is dissipative. By Lemma \ref{lem-w}, $M_\alpha
  \subseteq \dis$.
\end{enumerate}
In view of (\ref{d-alpha}) and the inclusion $T^{-1} \dis \subseteq \dis$, the
latter case implies $E_\alpha \subseteq \dis$, which is impossible
because $E_\alpha$ is a mixed \erg\ component by hypothesis. So the
statement about the \erg\ properties is proved.

Now consider $W := T^{-1} M_\alpha \setminus M_\alpha$. Certainly
$\leb(W)>0$, otherwise $T^{-1} M_\alpha = M_\alpha$ and $E_\alpha =
M_\alpha \subseteq \con$, which is false, again because $E_\alpha
\subseteq \icd$. This proves in particular that $\Z_\alpha$ is not an
isolated class, hence $\alpha \in \ater$. One readily checks that the
sets $\{ T^{-n} W \}_{n\in\N}$ are pairwise disjoint, so $W$ is
$T$-wandering. Moreover, $\bigcup_{n\in\N} T^{-n} W = E_\alpha
\setminus M_\alpha \subseteq \dis$.

For the first part of Proposition \ref{prop-ic-d} it remains to show
that $\omega(x) = M_\alpha$, for a.e.\ $x \in E_\alpha$.  But this
follows trivially from (\ref{d-alpha})-(\ref{t-alpha}) and the fact
that $T_\alpha$ is conservative and \erg.

\skippar

Viceversa, suppose $\alpha \in \ater$. This implies that $T^{-1}
M_\alpha \supset M_\alpha$ (strictly mod $\leb$, according to our
convention). The arguments used in the paragraph before the last 
one prove that
$\leb( E_\alpha \setminus M_\alpha ) > 0$ and $E_\alpha \setminus
M_\alpha \subseteq \dis$.

Now, set $A := M_\alpha \cap \Omega$. This is a $T_\alpha$-invariant
set, so the proof of Proposition \ref{prop-ic-c} applies, with
$T_\alpha$ in lieu of $T$. There are two cases:
\begin{enumerate}
\item $\leb(A)>0$. In this case, $A = M_\alpha$ is a conservative
  \erg\ component of $T_\alpha$---the unique component, in fact.  By
  Lemma \ref{lem-w}, $M_\alpha \subseteq \con$;
\item $\leb(A)=0$. In this case, a.e.\ point of $M_\alpha$ has an
  empty $\omega$-limit w.r.t.\ $T_\alpha$, equivalently, w.r.t.\ $T$.
  By Proposition \ref{prop-omega}, $M_\alpha \subseteq \id$.
\end{enumerate}
Observe that the two cases above correspond to the two cases described
in the first part of this proof. In any event, the first one gives
$E_\alpha \subseteq \icd$, and the second one gives $E_\alpha
\subseteq \id$.  As in the proof of Proposition \ref{prop-ic-c},
$\#\Z_\alpha < \infty$ implies the first case.  
\qed

\bigskip

\proofof{Theorem \ref{thm-dmho}} For $x \in \id = \R \setminus \Omega$
(cf.\ Proposition \ref{prop-omega}), one and only one of the following
occurs:
\begin{enumerate}
\item $\ds \lim_{n \to \infty} T^n(x) = +\infty$;
\item $\ds \lim_{n \to \infty} T^n(x) = -\infty$;
\item $\ds \limsup_{n \to \infty} T^n(x) = +\infty$ and $\ds 
  \liminf_{n \to \infty} T^n(x) = -\infty$.
\end{enumerate}
In the third case, (A5) implies that the \o\ of $x$ intersects
$[0,\rho]$ infinitely many times, thus having an accumulation point
there. This is a contradiction, and so $\id = \dis_{+\infty} \sqcup
\dis_{-\infty}$.

The remaining assertions will be proved only for $\dis_{+\infty}$, the
arguments for $\dis_{-\infty}$ being completely analogous.

Assume (A6) and $\leb(\dis_{+\infty}) > 0$. By definition of the
invariant set $\dis_{+\infty}$, the \o\ of a.e.\ $x \in
\dis_{+\infty}$ visits an infinite number of distinct intervals
$I_{j_n}$, with $\leb( I_{j_n} \cap \id ) > 0$. Here $(j_n)_n$ is a
subsequence of $\Z$ which depends on $x$. By (A6),
\begin{equation} \label{exa-135}
  \sum_n \leb(I_{j_n}) = \infty. 
\end{equation}
On the other hand, a.e.\ $x \in \dis_{+\infty}$ is also a density
point of $\dis_{+\infty}$. Hence, by Lemma \ref{lem-ex1} and
(\ref{exa-135}), $\leb(\dis_{+\infty}) = \infty$.

Lastly, we assume (A7) too and prove that $\dis_{+\infty}$ is an exact
component of $T$. We need the following distortion lemma.

\begin{lemma} \label{lem-ex2}
  Under the assumptions (A1)-(A7), there exists $D_1>1$ such that, for
  all measurable $A' \subseteq \R$ and all $j \in \Z$, $\leb(A' | I_j)
  > 1 - \delta$ implies $\leb( TA' | I_k ) > 1 - D_1 \delta$, for all
  $k_{1j} \le k \le k_{2j}$ (equivalently, for all $k$ such that $I_k
  \subset TI_j =: J_j$).
\end{lemma}

\proofof{Lemma \ref{lem-ex2}} Recall the meaning of the constants
$\theta, \theta_o, \bigj$; cf.\ (A1), (A5)-(A7) and following remarks.

For $A', j, k$ as in the statement of the lemma, set $B := I_j
\setminus A'$. By (A1), $\leb(B) < \delta \theta$.  By Corollary
\ref{cor-dist}, $T$ expands $B$ by a rate that is at most $D$ times
the average expansion rate of $I_j$:
\begin{equation} \label{exa-140}
  \frac{\leb(TB)} {\leb(B)} \le D \, \frac{\leb(TI_j)} {\leb(I_j)} \le
  D \, \frac{\bigj \theta} {\theta_o}.
\end{equation}
In the worst case, $TB$ lands entirely in $I_k$, whence $\leb(TB |
I_k) < D \bigj (\theta / \theta_o)^2 \delta$.  Setting $D_1 := D \bigj
(\theta / \theta_o)^2$ and noticing that $TA' \cap I_k \supseteq I_k
\setminus TB$ yields the desired result.  
\qed

\skippar

Back to the proof of Theorem \ref{thm-dmho}: given $A \subset \R$, we
say that a set of the type $C = \bigsqcup_{k=i}^{\ell} I_k$ is
\emph{$A$-prevalent} if $\ell-i+1 \ge \bigj$ and $\leb(A | I_k) >
1/2$, $\forall k \in \{ i, i+1, \ldots, \ell \}$. In other words, $C$
is made up of at least $\bigj$ Markov intervals, in each of which the
relative \me\ of $A$ is bigger than half; $\bigj$ is the positive
integer defined in (\ref{def-bigj}).

Suppose, by absurd, that $A, B \subset \dis_{+\infty}$ are disjoint,
$T$-invariant and of positive \me. We prove that, for a typical $x \in
A$, $x_n := T^n(x)$ belongs to an $A$-prevalent set, for all $n$ large
enough.

In fact, set $\delta := D_1^{-\bigj+1} /2$, where $D_1$ is the universal
constant provided by Lemma \ref{lem-ex2}. By Lemma \ref{lem-ex1} and
the invariance of $A$, $\exists n_1 =n_1(x)$ such that, $\forall n \ge
n_1$, $\leb(A | \ione{x_n}) > 1 - \delta$. Applying Lemma
\ref{lem-ex2} recursively $\bigj - 1$ times gives $\leb( A | I_k ) >
1/2$, for all $I_k \subset T^{\bigj-1} \ione{x_n}$. Observe that, by
(A2) and (A7), $T^{\bigj-1} \ione{x_n} = \bigsqcup_{k=i}^\ell I_k$,
for some $i < \ell$ that depend on $x$ and $n$. If we show that
\begin{equation} \label{exa-150}
  \ell-i+1 \ge \bigj, 
\end{equation}
we have proved that $T^{\bigj-1} \ione{x_n}$ is an $A$-prevalent set
that contains $x_{n+\bigj-1}$. This would give the assertion made in
the previous paragraph, because the above argument holds for
\emph{all} $n \ge n_1(x)$. But (\ref{exa-150}) is easily verified: by
(A2) and (A7), $( T^j \ione{x_n} )_{j \ge 0}$ is an increasing
sequence of sets that are unions of adjacent Markov intervals. The
sequence must be strictly increasing, otherwise, for some $j$, $T^j
\ione{x_n}$ would be forward-invariant, contradicting that $x_{n+j}
\to +\infty$, as $j \to \infty$. This implies that $T^{\bigj-1}
\ione{x_n}$ is made up of at least $\bigj$ intervals, which is
precisely (\ref{exa-150}).

So there are infinitely many $A$-prevalent sets in any right half-line
of $\R$.  On the other hand, applying the above to a typical $y \in
B$, we have that $y_n := T^n(y)$ belongs to a $B$-prevalent set, for
all large $n$. But $y_n \to +\infty$, and the distance between $y_n$
and $y_{n-1}$, in terms of intervals, is at most $\bigj-1$. This means
that, for $n$ big enough, a point $y_n$ must fall in an $A$-prevalent
set. But this is a contradiction, since $A$-prevalent sets and
$B$-prevalent sets cannot overlap.  Hence, $A$ and $B$ cannot be
disjoint and $\dis_{+\infty}$ is an \erg\ component.

To prove that it is also an exact component we apply Proposition
\ref{prop-mn} to $T|_{ \dis_{+\infty} }$. In fact, given $A \subseteq
\dis_{+\infty}$, the arguments used earlier show that, for a.e.\ $x
\in A$ and all large $n$, depending on $x$, $\leb(T^n A | \ione{x_n} )
> 1 - 1/2D_1 > 1/2$.  By Lemma \ref{lem-ex2} and (A7), $\leb(T^{n+1} A |
\ione{x_n} ) > 1/2$, whence $\leb(T^{n+1} A \cap T^n A ) >0$, which is
the hypothesis of Proposition \ref{prop-mn}.  
\qed

\subsection{Quasi-lifts and finite modifications}
\label{subs-pf-qlfm}

\proofof{Proposition \ref{prop-exa-ql}} Because of (A7) and
(\ref{trn-comm}), $T$ is aperiodic and Markov-indecomposable.  Recall
that for this map $\con = \ic$. If $T$ is not dissipative, then it is
conservative, irreducible and exact by Corollary
\ref{cor-c}\emph{(d)}.

If $T$ is dissipative, Proposition \ref{prop-no-do} and Theorem
\ref{thm-dmho} entail that there are at most two exact components:
$\dis_{+\infty}$ and $\dis_{-\infty}$. Suppose that neither has \me\
zero. By (\ref{trn-comm}), both components are invariant for the
action of $\trn$, thus $\leb (\dis_{\pm\infty} | I_j)$ is constant in
$j$.  But the proof of Theorem \ref{thm-dmho} shows that there exists
a (sufficiently large) $j$ such that $\leb (\dis_{+\infty} | I_j) >
1/2$. This, then, holds for all $j \in \Z$. The same can be proved for
$\leb (\dis_{-\infty} | I_j)$. It follows that our assumption was
wrong and there is only one exact component.

In order to characterize which type of exact component $\R$ one
obtains, depending on $\phi$, we look at $S_n \phi(y) :=
\sum_{k=0}^{n-1} \phi \circ T_a (y)$, the Birkhoff sum of $\phi$ for
the \dsy\ $( [0,a), \mu_a, T_a )$.  The collection of all these random
variables is also referred to as the (additive) \emph{cocycle}
generated by $\phi$. Any such cocycle is called \emph{recurrent} if,
for $\mu_a$-a.e.\ $y$,
\begin{equation} \label{def-co-rec}
  \liminf_{n \to \infty} |S_n \phi (y)| = 0.
\end{equation}   
A classical result by Atkinson \cite{at} shows that, if $T_a$ is \erg\
and $\phi$ is integrable w.r.t $\mu_a$---both holding here---then
$(S_n \phi)_{n \in \N}$ is recurrent \iff $\E_\mu (\phi) = 0$.

Suppose this is the case. The iterates of the $\Z$-extension $T_\phi$,
cf.\ Section \ref{subs-ql}, are of the form
\begin{equation} \label{pql-10}
  T_\phi^n (y, j) = (T_a^n (y), j + S_n \phi (y)).
\end{equation}
By (\ref{def-co-rec})-(\ref{pql-10}), the $T_\phi$-\o\ of $\leb$-a.e.\
point of $[0,a) \times \{ j \}$ has infinitely many returns
there. Passing to its conjugated map $T$, this means that a.e.\ $x \in
I_j$ has infinitely many returns to $I_j$. Thus, $T$ is
conservative. (This is obvious by the invariance of $\mu$, but would
hold anyway by the Markov properties of $T$, cf.\ proof of Proposition
\ref{prop-omega}.)

Suppose instead $\E_\mu (\phi) > 0$. The \erg ity of $T_a$ gives that,
for a.e.\ $y \in [0,a)$, $\lim_{n \to \infty} S_n \phi (y) =
+\infty$. Once again, (\ref{pql-10}) and the correspondance between
$T_\phi$ and $T$ prove that $\lim_{n \to \infty} T^n(x) = +\infty$ for
a.e.\ $x \in \R$, namely, $\R = \dis_{+\infty}$. Analogously for the
third case.  
\qed

\bigskip

\proofof{Proposition \ref{prop-exa-fm}} $T$ is Markov-indecomposable
by hypothesis and aperiodic by (A7). If $\con = \ic$ has positive \me,
then, by Corollary \ref{cor-c}\emph{(d)}, $T$ is conservative,
irreducible and exact, ending the proof of the Proposition.

Hence, let us assume that $\dis = \id = \R$. Using the notation of
Section \ref{subs-fm}, let us see what implications this has on the
dynamics of $T_o$.

We denote by $\con_o$ the conservative part of $T_o$, and by
$\dis^o_{\pm \infty}$ the sets defined in Theorem \ref{thm-dmho},
relative to $T_o$. By Proposition \ref{prop-exa-ql}, one of these sets
is the whole $\R$ and $T_o$ is exact. If $\R = \con_o$, \erg ity and
conservativity entail that the forward $T_o$-\o\ of a.e.\ $x \in \R$
intersects $B := \bigsqcup_{j=-k_o}^{k_o} I_j$. This occurs in
particular for a.e.\ $x \in \R \setminus B$, where $T_o$ and $T$
coincide, implying that the forward $T$-\o s of a.a.\ $x \in \R$
accumulate in $B$. This conclusion contradicts the dissipativity of
$T$, i.e., $\R = \id$.

Therefore, either $\dis^o_{+\infty}$ or $\dis^o_{-\infty}$ has full
\me. Suppose, w.l.g., that it is $\dis^o_{+\infty}$. We want to prove
that the same occurs for $\dis_{+\infty}$, equivalently, $\leb(
\dis_{-\infty} ) = 0$. Assume the contrary and define, for $\ell \in
\N$,
\begin{equation} 
  A_\ell:= \rset{x \in \dis_{-\infty}} {T^n(x) \in 
  \bigsqcup_{j < -k_o} I_j, \ \forall n \ge \ell }.
\end{equation}
Clearly, $A_\ell \subseteq A_{\ell+1}$ and $\bigcup_{\ell \in \N}
A_\ell= \dis_{-\infty}$. Thus, $\exists \ell$ such that $\leb( A_\ell)
> 0$.  On the other hand, $T^\ell A_\ell \subseteq \dis^o_{-\infty}$,
because $T^\ell A_\ell \subseteq T^\ell \dis_{-\infty}=
\dis_{-\infty}$ and $T^n(x) = T_o^n(x)$, $\forall n \ge \ell$ (as all
such points lie outside of $B$). The non-singularity of $T$ gives
$\leb(\dis^o_{-\infty}) \ge \leb(T^\ell A_\ell) > 0$, which is a
contradiction. Thus, our assumption was wrong and $\dis_{+\infty} =
\R$ mod $\leb$.  Finally, $T$ is exact by Theorem \ref{thm-dmho}.
\qed

As promised in Section \ref{subs-fm}, we present here a stronger
version of Proposition \ref{prop-exa-fm}, together with its proof.

\begin{proposition} \label{prop-exa-fm2} 
  Let $T$ be a finite modification of the quasi-lift $T_o$. Denote by
  $\mu_o$ the $\trn$-invariant \me\ preserved by $T_o$, and by
  $\phi_o$ discrete displacement \fn\ for $T_o$, as introduced in
  Section \ref{subs-ql}. Suppose that $T$ verifies (A1)-(A7) and
  preserves a Lebesgue-equivalent \me\ $\mu$ with the following
  properties:
  \begin{itemize}
  \item[(i)] $\exists \theta_2 > 0$ such that $\mu(I_j) \le \theta_2$,
    $\forall j \in \Z$;
  \item[(ii)] $\avg = \avg_o$ (Definition \ref{def-avg-avgo}).
  \end{itemize}
  Then $\R$ equals $\con$, $\dis_{+\infty}$, or $\dis_{-\infty}$,
  depending on $\E_{\mu_o} (\phi_o)$, the drift of $T_o$, being,
  respectively, zero, positive, or negative.
\end{proposition}

Observe that, unlike Proposition \ref{prop-exa-ql}, the case $\R =
\con$ does not guarantee that $T$ is exact (cf.\ Countexample 3 of
Appendix \ref{subs-ce}). But, if $T$ is also Markov-indecomposable,
exactness holds by Corollary \ref{cor-c}\emph{(d)}. This shows how
Proposition \ref{prop-exa-fm} is a corollary of Proposition
\ref{prop-exa-fm2}.

\bigskip

\proofof{Proposition \ref{prop-exa-fm2}} We use notation and several
arguments from the proof of Proposition \ref{prop-exa-fm}.

If $\E_{\mu_o} (\phi_o) = 0$, that is, $\con_o = \R$, the part of the
previous proof that shows that a.a.\ \o s accumulate in $B$ still
holds. Since $\icd$ is null by the invariance of $\mu$, it must be
$\con = \R$.

If $\E_{\mu_o} (\phi_o) > 0$, namely, $\dis^o_{+\infty} = \R$, the
argument given earlier whereby $\leb( \dis_{-\infty} ) = 0$ continues
to work. But $\leb( \dis_{+\infty} ) > 0$, because $\dis_{+\infty}$
coincides with $\dis^o_{+\infty} = \R$ on a large set on the ``right
end'' of $\R$. In order to prove that $\dis_{+\infty}$ has full \me,
we need to verify that $\leb(\con)=0$.

Suppose instead that $\leb(\con) > 0$. We show that $\exists k_1 \in
\Z$ such that
\begin{equation} \label{exa-fm2-5}
  \con = \bigsqcup_{j < k_1} I_j, \quad \dis_{+\infty} = 
  \bigsqcup_{j \ge k_1} I_j. 
\end{equation}
Recall that $\con = \ic$ is $\mar$-measurable by Corollary
\ref{cor-c}\emph{(a)}. If (\ref{exa-fm2-5}) does not hold, there exist
$j_2 \in \Z$ with $I_{j_2} \subset \con$, and a positive-\me\ set of
$x \in \dis_{+\infty}$ such that $x$ lies to the left of $I_{j_2}$ and
$T(x)$ lies to the right of $I_{j_2}$. Therefore, $\exists j_1 < j_2$
such that a positive-\me\ subset of such $x$ belong in
$I_{j_1}$. Since $\dis_{+\infty} = \R \setminus \con$ is also
$\mar$-measurable, $I_{j_1} \subseteq \dis_{+\infty}$, whence $T
I_{j_1} \subseteq \dis_{+\infty}$. But, for any $x \in I_{j_1}$ with
the properties stated earlier, (A7) shows that both $x$ and $T(x)$
belong in $T I_{j_1}$. Since $T I_{j_1}$ is an interval, it must
include $I_{j_2}$ too, which is absurd, because $I_{j_2} \subset
\con$. Thus (\ref{exa-fm2-5}) is established.

For all $\ell > k_1$, set $B_\ell := \bigsqcup_{j = k_1}^{\ell-1}
I_j$. The invariance of $\mu$ gives
\begin{equation} \label{exa-fm2-10} 
  \mu(B_\ell \setminus T^{-1} B_\ell) = \mu(T^{-1} B_\ell \setminus 
  B_\ell).
\end{equation}
The above l.h.s.\ comprises all the points that leave $B_\ell$ in one
iteration of $T$; the r.h.s.\ comprises all the points that enter
$B_\ell$ in one iteration of $T$. Now, recall the meaning of $k_{1,0},
k_{2,0}$ from (A2), and choose a sufficiently large $k_2 > k_1$ so
that $T$ and $T_o$ coincide on $\bigsqcup_{j \ge k_2 - k_{2,0}} I_j$.
From now on, we restrict to $\ell \ge k_2$.  Since $\bigsqcup_{j \ge
k_1} I_j = \dis_{+\infty}$ is $T$-invariant, points can only leave
or enter $B_\ell$ ``through its right end''. In formula:
\begin{align}
  B_\ell \setminus T^{-1} B_\ell \ &= \: \bigsqcup_{i=1}^{k_{2,0}} \
  \bigsqcup_{j=0}^{k_{2,0} - i} I_{\ell - i} \cap T^{-1} I_{\ell + j} ;
  \label{exa-fm2-20} \\
  T^{-1} B_\ell \setminus B_\ell \ &= \!\!\bigsqcup_{i=0}^{|k_{1,0}| - 1}
  \bigsqcup_{j=1}^{|k_{1,0}| - i} I_{\ell + i} \cap T^{-1} I_{\ell - j} .
  \label{exa-fm2-30}
\end{align}
We can replace $T$ with $T_o$ in the above r.h.sides, as already
observed. As $\ell$ varies, the resulting intervals are translation of
each other and it is relatively straightforward to evaluate their
$\mu_o$-\me s: translate the r.h.sides of (\ref{exa-fm2-20}) and
(\ref{exa-fm2-30}) via the maps $\trn^{-\ell+i}$ and $\trn^{-\ell-i}$,
respectively. This yields:
\begin{align}
  \mu_o ( B_\ell \setminus T^{-1} B_\ell ) &= \sum_{k=1}^{k_{2,0}} 
  k \, \mu_o ( I_0 \cap T_o^{-1} I_k ) ;
  \label{exa-fm2-40} \\
  \mu_o ( T^{-1} B_\ell \setminus B_\ell ) &= \sum_{k=1}^{|k_{1,0}|} 
  k \, \mu_o ( I_0 \cap T_o^{-1} I_{-k} ) .
  \label{exa-fm2-50}
\end{align}
Therefore, for all $\ell \ge k_2$,
\begin{equation} \label{exa-fm2-60}
\begin{split}
  \mu_o ( B_\ell \setminus T^{-1} B_\ell ) -
  \mu_o ( T^{-1} B_\ell \setminus B_\ell ) \: &=
  \sum_{k=-k_{1,0}}^{k_{2,0}} \!\!\!\ k \, \mu_o ( I_0 \cap T_o^{-1} 
  I_k ) \\ 
  &= \mu_o (I_0) \, \E_{\mu_o} (\phi_o) > 0,
\end{split}  
\end{equation}
as we have assumed.

Now consider the \fn\
\begin{equation}
  F := \sum_{j \in \Z} \left( 1_{B_{k_2} \setminus T^{-1} B_{k_2}} 
  - 1_{T^{-1} B_{k_2} \setminus B_{k_2}} \right) \circ \trn^j. 
\end{equation}
Once again, observe that there is no harm in replacing $T$ with $T_o$
in the above definition. Since $F$ is bounded and periodic, it is a
global \ob.  For $k \ge k_2$, set $V_k := \bigsqcup_{\ell = k_2}^k
I_\ell = [a k_2, a(k+1)]$. It is easily verified that
\begin{equation} \label{exa-fm2-65}
  F 1_{V_k} - \sum_{\ell=k_2}^k  \left( 1_{B_\ell \setminus 
  T^{-1} B_\ell} - 1_{T^{-1} B_\ell \setminus B_\ell} \right) 
\end{equation}
is a compactly supported bounded \fn\ whose $L^\infty$-norm is
constant in $k$ and whose support, though varying with $k$, is always
contained in at most $2( |k_{1,0}| + k_{2,0}) $ Markov
intervals. Therefore, by the hypothesis \emph{(i)}, the integral of
(\ref{exa-fm2-65}) is uniformly bounded in $k$. But it is a
requirement of Definition \ref{def-avg-avgo} that $\mu$ be an infinite
\me, whence
\begin{equation} \label{exa-fm2-70}
  \lim_{k \to +\infty} \mu_{V_k} (F) = \lim_{k \to +\infty} \frac1 
  {\mu(V_k)} \sum_{\ell=k_2}^k \left( \mu(B_\ell \setminus T^{-1} 
  B_\ell) - \mu(T^{-1} B_\ell \setminus B_\ell) \right) = 0,
\end{equation}
the last equality coming from (\ref{exa-fm2-10}). On the other hand,
via the periodicity of $\mu_o$ and $F$, and using (\ref{exa-fm2-60}),
\begin{equation} \label{exa-fm2-80}
  \lim_{k \to +\infty} (\mu_o)_{V_k} (F) = 
  \frac1 {\mu_o (I_0)} \left( 
  \mu_o (B_{k_2} \setminus T^{-1} B_{k_2}) - \mu_o (T^{-1} B_{k_2}
  \setminus B_{k_2} \right) > 0.
\end{equation}

The hypothesis \emph{(ii)} implies that the l.h.sides of
(\ref{exa-fm2-70}) and (\ref{exa-fm2-80}) are the same, which is a
contradiction. Therefore, the supposition $\leb(\con) > 0$ was wrong,
proving that $\dis_{+\infty} = \R$.

Analogously, $\E_{\mu_o} (\phi_o) < 0$ gives $\dis_{-\infty} = \R$.
\qed

\subsection{Infinite mixing}

\proofof{Theorem \ref{thm-mix-ql}} Assertion \emph{(a)} comes from
Proposition \ref{prop-glm1}, because $T$ is exact by Proposition
\ref{prop-exa-ql}.

As for \emph{(b)}, fix $F \in \go_1$. Using definition (\ref{def-ak}),
one verifies that, for all $k \in \Z^+$,
\begin{equation} \label{ql-10}
  \frac{ \mu((F \circ T^n) 1_{[0,ak]}) } {\mu( [0,ak] )} = \int_{I_0} 
  (\mathcal{A}_k F) \circ T^n \, d\mu_{I_0}.
\end{equation}
This can be seen by writing $1_{[0,ak]} = \sum_{j=0}^{k-1} 1_{I_0}
\circ \trn^{-j}$ and using the commutativity of $T$ and $\trn$, which
are both $\mu$-invariant. By definition of $\go_1$, for every
$\eps>0$, there exists a large enough $k$ such that, for all $n \in
\N$,
\begin{equation} \label{ql-20}
  \left| \int_{I_0} (\mathcal{A}_k F) \circ T^n \, d\mu_{I_0} 
  - \int_{I_0} F_a \circ T^n \, d\mu_{I_0} \right| \le \eps.
\end{equation}
Recalling from Section \ref{subs-ql} the definition of the
\me-preserving \dsy\ $(\sa, \mu_a, T_a)$, we see that, since $F_a$ is
$a$-periodic,
\begin{align}
  \int_{I_0} F_a \circ T^n \, d\mu_{I_0} &=
  \frac1 {\mu(I_0)} \int_{\sa} F_a \circ T_a^n \, d\mu_a \nonumber \\
  &= \frac1 {\mu(I_0)} \int_{\sa} F_a \, d\mu_a \nonumber \\
  &= \int_{I_0} F_a \, d\mu_{I_0} = \avg(F), \label{ql-30}
\end{align}
as per definition (\ref{im-110}) (with the slight abuse of notation
whereby the projection of $F_a$ to $\sa \cong I_0$ is still called
$F_a$).

The following lemma is an easy consequence of the exactness of $T$.

\begin{lemma} \label{lem-pmu36}
  If, for some $b \in \C$ and $\eps \ge 0$, the limit
  \begin{displaymath}
    \limsup_{n \to \infty} \left| \frac{ \mu((F \circ T^n) g) } {\mu(g)} 
    - b \right| \le \eps
  \end{displaymath}
  holds for some $g \in L^1$, with $\mu(g) \ne 0$, then it holds for
  all $g \in L^1$, with $\mu(g) \ne 0$.
\end{lemma}

\proof See \cite[Lem.~3.6]{lpmu}.
\bigskip

Equations (\ref{ql-10})-(\ref{ql-30}) show that the lemma can be
applied, for any $\eps>0$, with $b = \avg(F)$ and $g = 1_{[0,ak]}$,
with $k$ depending on $\eps$.  Therefore, for all $g$ with $\mu(g) \ne
0$,
\begin{equation}
  \lim_{n \to \infty} \left| \mu((F \circ T^n) g) - \avg(F) \mu(g) 
  \right| = 0.
\end{equation}
The case $\mu(g) = 0$ was already covered when we proved
\textbf{(GLM1)}. This ends the proof of \emph{(b)}.

By simple density arguments, \textbf{(GGM2)} will be verified once the
limit (GGM2) is proved for any pair of quasiperiodic \ob s $F,G$. More
precisely, assume $F \circ \trn = e^{\i a \beta}F$ and $G \circ \trn =
e^{\i a \gamma}G$. Notice that, if $\beta \ne 0$ mod $2\pi / a$, then
$\avg(F) = 0$. The analogous implication holds for $G$.

Set $g = G 1_{I_0} / \mu(I_0)$. Since $V = [ak, a(\ell+1)] =
\bigsqcup_{j=k}^\ell I_j$, we can write
\begin{equation} \label{ql-40}
  G 1_V = \sum_{j=k}^\ell G 1_{I_j} = \sum_{j=k}^\ell ((G \circ \trn^j) 
  1_{I_0}) \circ \trn^{-j} = \mu(I_0) \sum_{j=k}^\ell e^{\i a \gamma j} g 
  \circ \trn^{-j},
\end{equation}
whence, using the quasi-periodicity of $F$, the commutativity of $T$
and $\trn$, and the $\trn$-invariance of $\mu$,
\begin{align} \label{ql-50}
  \mu_V ((F \circ T^n)G) &= \frac1 {(\ell-k+1) \mu(I_0)} \int_\R 
  (F \circ T^n) G 1_V \, d\mu \nonumber \\
  &= \frac1 {\ell-k+1} \sum_{j=k}^\ell \int_\R (F \circ T^n) 
  e^{\i a \gamma j} (g \circ \trn^{-j}) \, d\mu \nonumber \\
  &= \frac{ \sum_{j=k}^\ell e^{\i a(\beta + \gamma)j} } {\ell-k+1}  
  \int_\R (F \circ T^n) g \, d\mu .
\end{align}
In the last term above, the factor in front of the integral is bounded
by 1 uniformly in $k, \ell$, namely, in $V$, while the integral does
not depend on it.  In fact, the latter term is $\mu((F \circ T^n) g)$
and, by \textbf{(GLM2)}, converges to $\avg(F) \mu(g)$, as $n \to
\infty$ (since $\go_2 \subset \go_1$).

We now have three cases:
\begin{enumerate}
\item \underline{$\beta \ne 0$ mod $2\pi / a$}. In this case $\avg(F)
  = 0$, therefore (\ref{ql-50}) converges to 0, as $n \to \infty$,
  uniformly in $V$. In particular, it converges to 0 in the joint
  infinite-volume and time limit; cf.\ (\ref{im-60}).
  
\item \underline{$\beta = 0$ mod $2\pi / a$ and $\gamma \ne 0$ mod
  $2\pi / a$}.  In this case $\avg(G) = 0$ and the factor in front
  of the integral in (\ref{ql-50}) vanishes when $\ell-k \to \infty$;
  that is, uniformly in $V$ as $\mu(V) \to \infty$, i.e., in the
  infinite-volume limit. On the other hand, the integral is bounded by
  $\| F \|_\infty \| g \|_1$, uniformly in $n$. This implies that,
  again, in the joint infinite-volume and time limit, (\ref{ql-50})
  converges to 0.

\item \underline{Both $\beta$ and $\gamma$ are 0 mod $2\pi / a$}. In
  this case the factor in front of the integral is identically 1,
  (\ref{ql-50}) no longer depends on $V$ and, for $n \to \infty$,
  tends to $\avg(F) \mu(g) = \avg(F) \avg(G)$, which is the same as
  the joint infinite-volume and time limit, here.
\end{enumerate}
In all these cases, the limit (GGM2) is verified.

In view of Proposition \ref{prop-mgg2-1}, \textbf{(GGM1)} will be
shown once we have proved that $\avg((F \circ T^n)G)$ exists for all
$F,G \in \go_2$ and $n \in \N$. Once again, since $\avg$ is a
continuous \fn al in the $L^\infty$-norm, it is enough to prove the
assertion for $F,G$ quasiperiodic. But, in that case, $F \circ T^n$ is
quasiperiodic by (\ref{trn-comm}), which implies the same for $(F
\circ T^n)G$, which thus has an infinite-volume average.  This
completes the proof of assertion \emph{(c)} of Theorem
\ref{thm-mix-ql}.  
\qed

\bigskip

\proofof{Proposition \ref{prop-mix-fm}} As in the previous proof,
assertion \emph{(a)} comes from Propositions \ref{prop-exa-fm} and
\ref{prop-glm1}.

Now for \emph{(b)}. Consider $F, G \in \go_2$ and fix $n \in \N$.
Since $T^n$ is a finite modification of $T_o^n$, as emphasized in
Section \ref{subs-fm}, the \fn\ $(F \circ T^n) G$ differs from $(F
\circ T_o^n) G$ by a compactly supported and bounded \fn\ of
$\R$. This shows that $\avg_o ( (F \circ T^n) G )$ exists \iff $\avg_o
( (F \circ T_o^n) G )$ does, and they are equal. Theorem
\ref{thm-mix-ql}\emph{(c)} then implies that
\begin{equation} \label{ql-60}
  \lim_{n \to \infty} \avg_o ( (F \circ T^n) G ) = \avg_o (F) \avg_o (G),
\end{equation}
for all $F,G \in \go_2$. By the hypothesis on $\mu$, the above holds
as well with $\avg$ in place of $\avg_o$, ending the proof of
\emph{(b)}.
\qed

\subsection{The example of the random walk}

\proofof{Proposition \ref{prop-fmrw-mix}} The irreducibility of $T$ is
apparent from the expression of $\qrw$, see (\ref{fmrw-70}) and
(\ref{rw-50}). The exactness then comes from Proposition
\ref{prop-exa-fm}. Moreover, the proof of Proposition
\ref{prop-exa-fm} (in Section \ref{subs-pf-qlfm}) shows that, under
its hypotheses, a finite modification of a conservative map is
conservative. In the case at hand, $T_o$ is clearly conservative,
whereby $T$ is as well.

Thus, Proposition \ref{prop-fmrw-mix} will be proved once the limit
(GLM2) is proved for any $F \in \go'$ and some $g \in L^1$, with
$\leb(g) \ne 0$; see Lemma \ref{lem-pmu36} \emph{et seq}.  We take $g$
to be of the form
\begin{equation} \label{fmrw-80}
  g_\pi := \sum_{j \in \Z} \pi_j 1_{I_j},
\end{equation}
where $\pi := (\pi_j)_{j \in \Z}$ is a symmetric (i.e., $\pi_j =
\pi_{-j}$, $\forall j \in \Z$) and half-monotonic (i.e., $\pi_j \ge
\pi_{j+1}$, $\forall j \in \N$) stochastic vector on $\Z$. In
particular, $g_\pi$ is a density w.r.t.\ the Lebesgue \me, namely,
$g_\pi \ge 0$ and $\leb(g_\pi)=1$.

The following lemma (which extends \cite[Lem.~7]{bcll} to the
non-homogeneous random walk at hand) states in particular that the set
of densities thus constructed is closed under the action of the
dynamics, that is, under the action of the Perron-Frobenius operator
$P$ introduced in (\ref{pf-10})-(\ref{pf-30}).

\begin{lemma} \label{lem-fmrw-decr}
  If $\pi$ is a symmetric and half-monotonic stochastic vector on $\Z$
  and $g_\pi$ is its corresponding density on $\R$ via
  (\ref{fmrw-80}), then $P^n g_\pi = g_{\pi^{(n)}}$, where $\pi^{(n)}
  := \pi \qrw^n$ is the evolution at time $n$ of the initial state
  $\pi$ for the random walk described above.  Moreover, $\pi^{(n)}$ is
  symmetric and half-monotonic.
\end{lemma}

\proofof{Lemma \ref{lem-fmrw-decr}} We prove all the assertions for
$n=1$ and the lemma will follow by induction. For the sake of the
notation, let us denote $\pi' := \pi^{(1)} = \pi \qrw$.

By (\ref{pf-30}) and (\ref{fmrw-80}), for all $x \in (k,k+1)$, $(P
g_\pi)(x) = \sum_j \pi_j q_{jk} =: \pi'_k$. This means that $P g_\pi =
g_{\pi'}$.  Also, by the symmetry properties of $\pi$ and $\qrw$,
\begin{equation} \label{fmrw-p10}
  \pi'_{-k} = \sum_j \pi_j q_{j,-k} = \sum_j \pi_{-j} q_{j,-k} =
  \sum_j \pi_j q_{-j,-k} = \sum_j \pi_j q_{jk} = \pi'_k.
\end{equation}
Finally, using both the symmetry and the half-monotonicity of $\pi$,
\begin{align} 
  \pi'_0 - \pi'_1 &= [( \pi_{-2} + \pi_{-1} + 5\pi_0 + \pi_1 + \pi_2) -
  (\pi_0 + 5\pi_1 + 2\pi_2 + \pi_3)]/9 \nonumber \\
  &= ( 4\pi_0 - 3\pi_1 - \pi_3 ) /9 > 0; \\[6pt]  
  \pi'_1 - \pi'_2 &=  [( \pi_0 + 5\pi_1 + 2\pi_2 + \pi_3) -
  (\pi_0 + 2\pi_1 + 3\pi_2 + 2\pi_3 +\pi_4)]/9 \nonumber \\
  &= (3\pi_1 - \pi_2 - \pi_3 - \pi_4)/9 > 0;
\end{align}
and, for $k \ge 2$,
\begin{align} 
  \pi'_k - \pi'_{k+1} &= [( \pi_{k-2} + 2\pi_{k-1} + 3\pi_k +
  2\pi_{k+1} + \pi_{k+2}) \nonumber \\
  &\phantom{= [} - ( \pi_{k-1} + 2\pi_k + 3\pi_{k+1} + 2\pi_{k+2} +
  \pi_{k+3})]/9 \nonumber \\
  &= ( \pi_{k-2} + \pi_{k-1} + \pi_k - \pi_{k+1} - \pi_{k+2} -
  \pi_{k+3} )/9 >0.
\end{align}
Therefore $\pi'$ is decreasing on $\N$ and the proposition is proved.
\qed

\bigskip

Now for the core argument. Without loss of generality we assume that
$\lavg'(F) = 0$ (for (GLM2) is trivial when $F$ is a constant).  Set
$f_j := \int_j^{j+1} F \, d\leb$. The assumption implies that $\forall
\eps>0$, $\exists \ell \in \N$ such that, $\forall k \ge \ell$,
\begin{equation} \label{fmrw-p50}
  \frac1{2k+1} \left| \sum_{j=-k}^k f_j \right| = \left| 
  \leb_{[-k,k+1]}(F) \right| \le \frac{\eps}2.
\end{equation}
By (\ref{pf-10}) and Lemma \ref{lem-fmrw-decr} we have
\begin{equation} \label{fmrw-p60}
\begin{split}
  \int_\R (F \circ T^n) g_\pi \, d\leb &= \int_\R F g_{\pi^{(n)}} \, 
  d\leb \\
  &= \sum_{j \in \Z} f_j \, \pi_j^{(n)} \\
  &= \sum_{k=0}^\infty \left(\pi_k^{(n)} - \pi_{k+1}^{(n)} \right)
  \sum_{j=-k}^k f_j \\
  &=: S_\ell + S_\ell' ,
\end{split}
\end{equation}
where $S_\ell$ and $S_\ell'$ correspond to restricting the outer
summation to $\sum_{k=0}^{\ell-1}$ and $\sum_{k=\ell}^{\infty}$,
respectively. Observe that the third equality of (\ref{fmrw-p60})
comes from disintegrating the density $( \pi_j^{(n)} )_j$ in
``horizontal slices'' of width $2k+1$ and height $\pi_k^{(n)} -
\pi_{k+1}^{(n)}$.

By (\ref{fmrw-p50}) we obtain
\begin{equation} \label{fmrw-p70}
\begin{split}
  \left| S_\ell' \right| &\le \sum_{k=\ell}^\infty \left(\pi_k^{(n)} - 
  \pi_{k+1}^{(n)} \right) \left| \sum_{j=-k}^k f_j \right| \\
  &\le \frac{\eps}2 \sum_{k=\ell}^\infty \left(\pi_k^{(n)} - 
  \pi_{k+1}^{(n)} \right) (2k+1) \le \frac{\eps}2,
\end{split}
\end{equation}
because $\sum_{k \in \N} (\pi_k^{(n)} - \pi_{k+1}^{(n)}) (2k+1) =
\sum_{j \in \Z} \pi_j^{(n)} = 1$, as in (\ref{fmrw-p60}).

To estimate $S_\ell$ we need a property of the \dsy\ which is an easy
consequence of its exactness.

\begin{lemma} \label{lem-llm}
  For all $f \in L^\infty \cap L^1$ and $g \in L^1$,
  \begin{displaymath}
    \lim_{n \to \infty} \leb( (f \circ T^n)g ) = \lim_{n \to \infty} 
    \leb( f (P^n g) ) = 0.
  \end{displaymath}
\end{lemma}

\proof See \cite[Thm.~3.5\emph{(b)}]{lpmu}. (In our terminology, see
\cite{lpmu}, the above notion is called \emph{local-local \m}, or
\textbf{(LLM)}, and is easily seen to be equivalent to the zero-type
property of \cite{hk}.)

\bigskip

Let us apply Lemma \ref{lem-llm} with $f = 1_{[-\ell+1, \ell]}$ and $g
= g_\pi$. There exists $N \in \N$ such that, $\forall n \ge N$,
\begin{equation} \label{fmrw-p80}
  \sum_{j=-\ell+1}^{\ell-1} \pi_j^{(n)} = \int_{-\ell+1}^\ell \!
  g_{\pi^{(n)}} \, d\leb = \int_{-\ell+1}^\ell \!\! P^n g_\pi \, d\leb \le 
  \frac{\eps} {2 \| F \|_\infty}.
\end{equation}
Therefore
\begin{equation} \label{fmrw-p90}
\begin{split}
  \left| S_\ell \right| &\le \sum_{k=0}^{\ell-1} \left(\pi_k^{(n)} -
  \pi_{k+1}^{(n)} \right) \left| \sum_{j=-k}^k f_j \right| \\
  &\le \| F \|_\infty \sum_{k=0}^{\ell-1} \left(\pi_k^{(n)} -
  \pi_{k+1}^{(n)} \right) (2k+1) \\
  &\le \| F \|_\infty \! \sum_{j=-\ell+1}^{\ell-1} \pi_j^{(n)} \le
  \frac{\eps}2 .
\end{split}
\end{equation}
Since $N$ is chosen depending on $\ell$, which in turns depends on
$\eps$, (\ref{fmrw-p60}), (\ref{fmrw-p70}) and (\ref{fmrw-p90}) prove
the assertion.  
\qed

\appendix

\section{Appendix}

\subsection{The importance of certain assumptions}
\label{subs-ce}

In this section we present some examples---or rather
counterexamples---of maps which clarify the role of some of our less
obvious assumptions. We refer in particular to Theorem \ref{thm-dmho},
which describes the exact components of $\id$, and Propositions
\ref{prop-exa-fm} and \ref{prop-exa-fm2}, which apply to finite
modifications of quasi-lifts. All the maps we present are Markov maps
associated to random walks.

\paragraph{Counterexample 1:} With reference to Theorem
\ref{thm-dmho}, this example illustrates the relevance of (A7) for the
exactness properties of $T$ on $\id$.

Let $T$ be the map associated to the random walk with the following
transition probabilities:
\begin{equation}
\begin{array}{rl}
  \forall j \le -1, & q_{j,j-2} = 1/2, \ q_{j,j} = q_{j,j+2} = 1/4; \\
  \forall j \in \{0,1\}, & q_{j,-2} = q_{j,-1} = q_{j,0} = q_{j,1} = 
  q_{j,2} = q_{j,3} = 1/6; \\
  \forall j \ge 2, & q_{j,j-2} = q_{j,j} = 1/4, \ q_{j,j+2} = 1/2.
\end{array}
\end{equation}
All other $q_{jk}$ are necessarily null. As seen in Section
\ref{subs-rw}, $T$ verifies (A1)-(A6). It is also irreducible.

Denote $\R_\mathrm{even} := \bigsqcup_{j \in \Z} [2j, 2j+1)$ and
$\R_\mathrm{odd} := \R \setminus \R_\mathrm{even}$. For $\epsilon \in
\{ \mbox{even, odd} \}$ indicate with $\dis_{\pm \infty}^{\epsilon}$
the set of all $x \in \R$ such that $\lim_{n \to \infty} T^n(x) \to
\pm \infty$ and $T^n(x) \in \R_\epsilon$ for all sufficiently large
$n$. This defines four invariant sets. Once we prove that all of them
have infinite Lebesgue \me, we have shown that $\id$ has at least four
\erg\ components. This demonstrates that, if (A1)-(A6) hold but (A7)
does not, the last assertion of Theorem \ref{thm-dmho} fails, even for
an irreducible $T$.

So let us verify that $\leb( \dis^\mathrm{even}_{+\infty} ) = \infty$,
the proof for the other sets being analogous. Consider the map $T_1$
corresponding to the homogeneous random walk $q_{j,j-2} = q_{j,j} =
1/4$, $q_{j,j+2} = 1/2$, $\forall j \in \Z$.  Clearly $T_1^{-1}
\R_\mathrm{even} = \R_\mathrm{even}$ and, by the transience of the
random walk, $\lim_{n \to \infty} T_1^n(x) \to +\infty$ for all
$x$. This implies that there is an infinite-\me\ set of $x \in
\R_\mathrm{even}$ such that $T_1^n(x) \ge 2$, $\forall n \in
\N$. Since $T$ and $T_1$ coincide on $\R_{\ge 2}$, all such $x$ belong
to $\dis^\mathrm{even}_{+\infty}$, which therefore has infinite
Lebesgue \me. (A more refined analysis, using the arguments of Section
\ref{subs-pf-ex}, would show that $\con = \icd = \emptyset$, and the
four sets $\{ \dis_{\pm \infty}^{\epsilon} \}$ are the exact
components of $T$.)

\paragraph{Counterexample 2:} The assumption (A5) is even more
important than (A7) in Theorem \ref{thm-dmho}. The following example
shows that if (A5) does not hold---even if all other assumptions
do---the tail $\sigma$-algebra of a Markov map can be very large,
making the \sy\ very far from being exact.

Consider the map $T$ representing the following random walk: 
\begin{equation}
\begin{array}{rl}
  \forall j \le 1, & q_{j,j} = q_{j,j+1} = 1/2 ; \\
  \forall j \ge 2, & q_{j,j} = q_{j,j+1} = \ldots = q_{j, 2^{j^2} -1} = 
  1/(2^{j^2} - j).
\end{array}
\end{equation}
In other words, 
\begin{equation} \label{rw-100}
  T(x) = \left\{
  \begin{array}{lll}
    j + 2 (x-j), & x \in [j, j+1), & j \le 1; \\
    j + (2^{j^2}-j) (x-j), & x \in [j, j+1), & j \ge 2.
  \end{array}
  \right.
\end{equation}
The intervals $I_j = [j, j+1]$ provide an alternative Markov partition
for $T$---see Remark \ref{rk-rw}---relative to which the \sy\ verifies
(A1)-(A4), (A6)-(A7). Evidently, it does not verify (A5).

Using the notation of (A2), let us denote by $\br_j$ the branch of $T$
over $I_j$, which is expressed by the r.h.s.\ of (\ref{rw-100}). For
$j \ge 3$, set
\begin{equation} \label{rw-110}
  X_j := \br_j^{-1} \left( [ 2^{(j-1)^2}, 2^{j^2} ) \right) = \left[
  j + \frac{2^{(j-1)^2} - j} {2^{j^2} - j} , j+1 \right) \subset I_j.
\end{equation}
Define also $X := \bigsqcup_{j \ge 3} X_j$ and $Y := \bigcap_{n \ge 0}
T^{-n} X$. By construction, $T|_X$ is a bijection $X \into \R_{\ge
16}$, hence $T|_{Y}$ is an invertible self-map of $Y$. See
Fig.~\ref{ce2}.

\fig{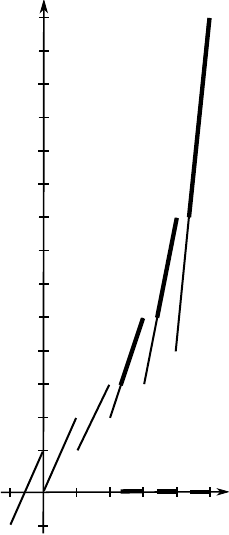}{3.7cm}{A rough sketch of the map $T$ of Counterexample
  2. The bold segments on the abscissa indicate the set $X$. The bold
  parts of the graph of $T$ represent $T|_X$, which is invertible.}

\begin{lemma} 
  $\leb(Y) > 0$. 
\end{lemma}

\proof It will suffice to show that $\leb(I_3 \cap Y) > 0$. Set $X_j'
= I_j \setminus X_j$ and $X' := \bigsqcup_{j \ge 3} X_j'$. By
(\ref{rw-110}),
\begin{equation} \label{rw-115}
  \leb(X_j') < \frac{2^{(j-1)^2} - j} 
  {2^{j^2} - j} <  \frac{2^{(j-1)^2}} {2^{j^2}} = 2^{-2j+1}.
\end{equation}

For every $n \ge 1$, $T^n$ acts as a piecewise linear bijection $I_3
\cap \bigcap_{i=0}^{n-1} T^{-i} X \into [k_n, \ell_n+1)$, for some
$k_n, \ell_n \ge 3$. Let us call this map $L_n$.  We have:
\begin{equation} \label{rw-120}
  I_3 \cap Y = X_3 \cap \bigcap_{n=1}^\infty T^{-n} X = X_3 \cap 
  \bigcap_{n=1}^\infty L_n^{-1} ( X \cap [k_n, \ell_n+1) ).
\end{equation}
The integers $k_n, \ell_n$ can be calculated recursively from the 
definition of $X_j$. For example, $k_1 = 2^4$, $k_2 = 2^9-1$, $k_2
= 2^{(2^4 - 1)^2}$, $k_2 = 2^{(2^9 - 1)^2}-1$, etc.  A very generous
lower bound for $k_n$ is $n$. In view of
(\ref{rw-115})-(\ref{rw-120}), the complementary set of $Y$,
w.r.t.\ $I_3$, measures
\begin{equation} \label{rw-130}
\begin{split}
  \leb\! \left( I_3 \cap \bigcup_{n=0}^\infty T^{-n} X' \right) 
  &\le \leb( X_3' ) + \sum_{n=1}^\infty \leb\! \left( L_n^{-1} 
  ( X' \cap [k_n, \ell_n+1) ) \right) \\
  &= \frac{13}{509} + \sum_{n=1}^\infty \sum_{j=k_n}^{\ell_n}
  \leb( L_n^{-1} X_j' ).
\end{split}
\end{equation}
On the other hand, 
\begin{equation}\label{rw-140}
\begin{split}
  \sum_{j=k_n}^{\ell_n} \leb( L_n^{-1} X_j' ) &= 
  \sum_{j=k_n}^{\ell_n} \leb( L_n^{-1} X_j' \,|\, L_n^{-1} I_j) \, 
  \leb( L_n^{-1} I_j ) \\
  &= \sum_{j=k_n}^{\ell_n} \leb( X_j' \,|\, I_j) \, 
  \leb( L_n^{-1} I_j ) \\
  &< \sum_{j=k_n}^{\ell_n} 2^{-2j+1} \, \leb( L_n^{-1} I_j ) \\
  &< 2^{-2n+1}.
\end{split}
\end{equation}
In the above we have used, from top to bottom: the affinity of
$L_n^{-1}|_{I_j}$, the inequality (\ref{rw-115}), the lower bound $n
\le k_n \le j$, and the fact that $L_n^{-1} [k_n, \ell_n+1) = I_3 \cap
\bigcap_{i=0}^{n-1} T^{-i} X$, whose Lebesgue \me\ is less than 1.

The estimate (\ref{rw-140}) shows that the l.h.s.\ of (\ref{rw-130})
is less than 1, which is equivalent to $\leb(I_3 \cap Y) > 0$, as
claimed. (Incidentally, by the same arguments as above, one proves
that $\leb(Y \cap I_\ell) \ge c$, for some $c>0$ and all $\ell \ge 3$,
showing that, in fact, $\leb(Y) = \infty$.)  
\qed

The bijectivity of $T|_Y : Y \into Y$ ensures that, for all $B
\subseteq Y$ and all $n \in \N$, $B = T^{-n} T^n B$, implying that $B
\in \sct$, the tail $\sigma$-algebra of $T$ defined in
(\ref{tail}). Equivalently, using the notation introduced in Section
\ref{subs-pf-ex}, $\sct \cap Y = \bor \cap Y$. In more suggestive
terms, $T$ cannot lose memory about $Y$, as it is invertible there!

\paragraph{Counterexample 3:} The next map is as much a counterexample
for Proposition \ref{prop-exa-fm} as an example for Proposition
\ref{prop-exa-fm2}. It is a finite modification of a quasi-lift which
preserves the same \me\ as the quasi-lift, but fails to be exact
because it is not Markov-indecomposable or, which is the same here,
because it has more than one conservative \erg\ component.

Let $T$ correspond to the random walk given by:
\begin{equation}
\begin{array}{rl}
  & q_{-1,-2} = 1/3, \ q_{-1,-1} = 2/3; \\  
  & q_{0,0} = 2/3, \ q_{0,1} = 1/3; \\
  \forall j \not\in \{ -1,0 \}, & q_{j,j-1} = q_{j,j} = q_{j,j+1} = 1/3.
\end{array}
\end{equation}
$T$ is a finite modification of a map $T_o$ which is associated to a
homogeneous random walk. The latter is thus a quasi-lift. By means of
Proposition \ref{prop-ds-im}, one readily checks that both $T$ and
$T_o$ preserve $\leb$. On the other hand, $T^{-1} \R^\pm = \R^\pm$,
showing that $T$ has at least two exact components. In fact, it is not
hard to see that $T$ is conservative and $E_1 := \R^+$, $E_2 := \R^-$
are the only two exact components. Since these sets are unions of
Markov intervals, $T$ is Markov-decomposable.

\paragraph{Counterexample 4:} Finally, we show that the preservation
of a \me\ equivalent to Lebesgue is also a crucial hypothesis for both
Propositions \ref{prop-exa-fm} and \ref{prop-exa-fm2}.

Let $T$ be given by this random walk:
\begin{equation}
\begin{array}{rl}
  & q_{-1,-2} = 1/3, \ q_{-1,-1} = 2/3; \\  
  & q_{1,1} = 2/3, \ q_{1,2} = 1/3; \\
  \forall j \not\in \{ -1,1 \}, & q_{j,j-1} = q_{j,j} = q_{j,j+1} = 1/3.
\end{array}
\end{equation}
This is a finite modification of the same $T_o$ as in the previous
case. In this case, however, $T$ does not preserve an
$\leb$-equivalent \me, as can be seen, e.g., by the fact that $[0,1]$
has only one inverse branch, which contracts by a factor 3.

What happens is that $T|_{\R^-}: \R^- \into \R^-$ and $T|_{\R_{\ge
1}}: \R_{\ge 1} \into \R_{\ge 1}$ are exact. Also, the Markov
interval $I_{0,-1} = [0, 1/3]$ feeds $\R^-$; $I_{0,1} = [2/3, 1]$
feeds $\R_{\ge 1}$; and $I_{0,0} = [1/3, 2/3]$ feeds both $\R^-$ and
$\R_{\ge 1}$. There are only two communicating classes, $\Z_1$ and
$\Z_2$, which are terminal. The corresponding sets are, respectively,
$M_1 = \R^-$ and $M_2 = \R_{\ge 1}$, with basins $E_1 = \R_{< 1/2}$
and $E_2 = \R_{> 1/2}$; cf.\ (\ref{m-alpha}),
(\ref{d-alpha})-(\ref{t-alpha}). Lastly, $T$ is Markov-indecomposable
because $I_{0,0}$ feeds both $M_1$ and $M_2$.

In conclusion, $T$ is not exact and $\R = \icd$, contrary to both the
statements of Proposition \ref{prop-exa-fm} and Proposition
\ref{prop-exa-fm2}.

\subsection{Distortion}
\label{subs-dist}

In this section of the Appendix we prove a standard distortion result
that is used in Section \ref{subs-pf-ex}.

\begin{lemma} \label{lem-dist}
  Under the hypotheses of Section \ref{sec-setup}, and using the
  notation (\ref{exa-5}), there exists $D>1$ such that, for all $n \ge
  1$, $\bj \in \Z^n$, $x,y \in I^{(n)}_\bj$,
  \begin{displaymath}
    D^{-1} \le \frac{ |(T^n)'(x)| } { |(T^n)'(y)| } \le D.
  \end{displaymath}
\end{lemma}

\proof For $0 \le k \le n$, set $x_k := T^k(x)$ and $y_k := T^k(y)$.
 
\skippar
\noindent
\textbf{Convention.}\ The notation $(x_k, y_k)$ denotes both the
interval $(x_k, y_k)$, when $x_k < y_k$, and the interval $(y_k, x_k)$,
when $y_k < x_k$.  
\skippar

By definition of $I^{(n)}_\bj$, see (\ref{exa-5}), we have
\begin{equation} \label{dist-5}
  (x_k, y_k) \subset I_{(j_k, \ldots, j_{n-1})}, 
\end{equation}
so that $T$ is twice differentiable on $(x_k, y_k)$. Therefore,
\begin{equation} \label{dist-10}
\begin{split}
  \log \frac{ |(T^n)'(x)| } { |(T^n)'(y)| } &= \sum_{k=0}^{n-1}
  \log \frac{ |T'(x_k)| } { |T'(y_k)| } \\
  &= \sum_{i=1}^n \left( \log |T'( \br_{j_{i-1}}^{-1} (x_i))| - 
  \log |T'( \br_{j_{i-1}}^{-1} (y_i))| \right) \\
  &= \sum_{i=1}^n \left| \frac{T''( \br_{j_{i-1}}^{-1} (z_i) )} 
  {T'( \br_{j_{i-1}}^{-1} (z_i) )} \, \frac1 {T'( \br_{j_{i-1}}^{-1} 
  (z_i) )} \right| (x_i - y_i),
\end{split}
\end{equation}
where $z_i \in (x_i, y_i)$. By (\ref{dist-5}), using (A1) and (A3), we
get $|x_i - y_i| \le \theta \lambda^{n-i}$, which, using (A4) as well,
gives
\begin{equation} \label{dist-20} 
  \left| \log \frac{ |(T^n)'(x)| } { |(T^n)'(y)| } \, \right| \le 
  \sum_{i=1}^n \eta \, \theta \, \lambda^{n-i} \le \frac{\eta \, \theta} 
  {1 - \lambda}.
\end{equation}
Renaming the rightmost term of (\ref{dist-20}) $\log D$ yields the
assertion.  
\qed

\begin{corollary} \label{cor-dist}
  Let $\bj = (j_0, j_1, \ldots, j_n) \in \Z^{n+1}$ be such that
  $\leb(I^{(n+1)}_\bj) > 0$. If $B \subset I^{(n+1)}_\bj$, then $T^n B
  \subset I_{j_n}$ and
  \begin{displaymath}
    \frac{ \leb(T^n B) } { \leb(I_{j_n}) } \le D \, \frac{ \leb(B) } 
    { \leb(I^{(n+1)}_\bj) }.
  \end{displaymath}
\end{corollary}

\proof By construction, $T^n I^{(n+1)}_\bj = I_{j_n}$, giving the
first assertion. As for the second,
\begin{equation} \label{dist-30}
  \frac{ \leb(T^n B) } { \leb(T^n I^{(n+1)}_\bj) } = \frac{ \int_B
  |(T^n)'(x)| \, dx } { \int_{I^{(n+1)}_\bj} |(T^n)'(x)| \, dx } \le
  \frac{ \max_B |(T^n)'| } { \min_{I^{(n+1)}_\bj} |(T^n)'| } \, \frac{
  \leb(B) } { \leb(I^{(n+1)}_\bj) },
\end{equation}
which, through Lemma \ref{lem-dist}, yields the corollary.  
\qed

\end{document}